\documentclass[11pt,twoside]{article}

\usepackage[
a4paper,
left=28.8mm,
right=28.8mm,
top=25mm,
bottom=25mm,
headheight=14pt,
headsep=16pt
]{geometry}

\usepackage{amsmath,amssymb,amsthm}
\usepackage{mathrsfs}

\usepackage{hyperref}
\usepackage[nameinlink,noabbrev]{cleveref}

\usepackage{fancyhdr}
\usepackage{cite}
\usepackage{enumitem}

\numberwithin{equation}{section}

\theoremstyle{plain}

\newtheorem{theorem}{Theorem}[section]
\newtheorem{corollary}[theorem]{Corollary}
\newtheorem{lemma}[theorem]{Lemma}

\theoremstyle{remark}

\theoremstyle{definition}

\newtheorem*{example}{Example}

\crefname{theorem}{Theorem}{Theorems}
\Crefname{theorem}{Theorem}{Theorems}

\crefname{lemma}{Lemma}{Lemmas}
\Crefname{lemma}{Lemma}{Lemmas}

\crefname{proposition}{Proposition}{Propositions}
\Crefname{proposition}{Proposition}{Propositions}

\crefname{corollary}{Corollary}{Corollaries}
\Crefname{corollary}{Corollary}{Corollaries}

\crefname{remark}{Remark}{Remarks}
\Crefname{remark}{Remark}{Remarks}

\crefname{definition}{Definition}{Definitions}
\Crefname{definition}{Definition}{Definitions}

\crefname{example}{Example}{Examples}
\Crefname{example}{Example}{Examples}

\crefname{equation}{equation}{equations}
\Crefname{equation}{Equation}{Equations}

\newcommand{\shorttitle}{Minkowski Problems for Mixed Centro-Section Integrals}
\newcommand{\shortauthors}{RENQING YOU}

\newcommand{\runningheadfont}{\fontsize{9}{11}\selectfont\scshape}

\title{Minkowski Problems for Mixed Centro-Section Integrals}

\author{
	Renqing You
}

\date{}

\date{}
\begin{document}
	\maketitle
	
	\begin{abstract}
	Recent studies in affine convex geometry have introduced two new families of geometric measures: affine dual curvature measures and centro-section measures. Motivated by these developments, we introduce mixed centro-section integrals and establish their variational formulas. Under the subspace concentration condition, we prove the existence of solutions to the corresponding Minkowski-type problems. We further introduce mixed bi-dual $m$-intersection bodies of index $q$ and establish the corresponding sharp affine isoperimetric inequalities. The resulting mixed centro-section measures include centro-section measures and affine dual curvature measures as special cases.
	\end{abstract}
	
	

	\bigskip
	\section{Introduction}
	\quad In recent years, Minkowski problems for geometric measures and the associated isoperimetric inequalities have attracted considerable attention in convex geometric analysis. Given a finite Borel measure on the unit sphere, a Minkowski-type problem asks whether the measure can be realized as a geometric measure of a convex body and, if so, how the
	body can be reconstructed. The classical Minkowski problem concerns the
	surface area measure of a convex body. Other fundamental examples include the
	Christoffel--Minkowski problem for mixed area measures
	\cite{2003-Guan-invent,2006-Guan-CPAM} and the Aleksandrov problem for integral curvature measures
	\cite{Aleksandrov1942,2019-SW-JEMS,Zhao-2019-PAMS}. For regular convex bodies, the corresponding
	Minkowski problems are equivalent to fully nonlinear Monge--Amp\`ere-type
	equations on the sphere; see, for example, the works of Nirenberg
	\cite{1953-Nirenberg-CPAM}, Cheng--Yau \cite{1976-CY-CPAM}, Pogorelov
	\cite{1978-Pogorelov-DC}, and Caffarelli
	\cite{1990-Caffarelli-Ann.mtah}. For a comprehensive account of the classical
	theory, we refer to Schneider's book \cite{2014-Schneider-CUP}, while a
	recent survey of Minkowski problems for geometric measures is given by Huang–Yang–Zhang \cite{2025-Huang-BAMS}. Further developments including the
	logarithmic Minkowski problem \cite{2002-SA-Advmath,2003-SA-Advmath,2013-Boroczky-JAMS,2019-CLZ-TAMS,2018-Zhu-advmath}, chord Minkowski problem \cite{2024-Lutwak-CPAM,2023-Xi-ANS}, $L_p$ Minkowski problem \cite{2006-CW-Advmtah,2018-Huang-AdvMath,1993-Lutwak-JDG,2018-Lutwak-Advmath, 2015-Zhu-JFA}, dual Minkowski problem\cite{2018-Huang-AdvMath,2016-Huang-ActaMath,2019-SW-JEMS,2017-ZY-CVPDE,2018-ZY-advmath}, Orlicz Minkowski problem 
	\cite{2019-Gardner-CVPDE,2020-Gardner-CVPDE} and centro-affine Minkowski problem \cite{2015-Zhu-JDG,2006-CW-Advmtah,2019-Li-IMRN}.
	
	Over the past three decades, the introduction of new geometric measures has substantially expanded the scope of the Brunn–Minkowski theory. Important
	examples include the cone-volume measure \cite{2013-Boroczky-JAMS}, the
	$L_p$ surface area measures \cite{1993-Lutwak-JDG}, the dual curvature
	measures \cite{2016-Huang-ActaMath}, and the chord measures
	\cite{2024-Lutwak-CPAM}. More recently, Cai, Leng, Wu, and Xi introduced the
	affine dual curvature measures \cite{2025-Cai-AdvMath} and the centro-section
	measures \cite{2026-Cai-AdvMath}. These developments illustrate the
	fundamental role of variational methods in geometric measures and associated Minkowski problems.
	
	The dual Brunn--Minkowski theory was introduced by Lutwak in the mid-1970s\cite{Lutwak-1975}. Many fundamental geometric invariants in the classical Brunn--Minkowski theory, including mixed volumes and quermassintegrals, together with their related formulas and geometric inequalities, have counterparts in the dual theory; see, e.g., \cite{2019-Boroczky-AdvMath, 2017-CLZ-JDE, 2006-CW-Advmtah, 2019-Gardner-CVPDE, 2018-HP-AdvMath, 2018-Huang-AdvMath, 2019-SW-JEMS, 1994-ZGY-TAMS, 2017-ZY-CVPDE}. In this dual setting, dual quermassintegrals play a particularly important role and are naturally expressed in terms of volumes of central sections.
	
	In \cite{2016-Huang-ActaMath}, Huang-Lutwak-Yang-Zhang established the celebrated variational formula for the dual quermassintegrals, and attacked the associated dual Minkowski problems. The recent breakthrough of Milman and Yehudayoff \cite{2023-Milman-JAMS},
	who established the first family of affine isoperimetric inequalities
	conjectured by Lutwak in 1988 \cite{1988-Lutwak-Math.Ann} for his affine
	quermassintegrals, has revived interest in these functionals. Inspired by these two works, Cai, Leng, Wu, and Xi \cite{2026-Cai-AdvMath} considered the $q$-th moments of the volumes of the intersections of a convex body with $m$-dimensional subspaces. These functionals are called the $m$-th centro-section integrals of
	index $q\in\mathbb{R}$ and are denoted by $\widetilde{\Psi}_{m,q}$. For $m \in \{1,\dots,n-1\}$, define
	\begin{equation}
		\widetilde{\Psi}_{m,q}(K) = \int_{G(n,m)} \operatorname{vol}_m(K \cap \xi)^q d\nu_m(\xi), \label{1-1}
	\end{equation}
	where $G(n,m)$ denotes the Grassmannian manifold of $m$-dimensional subspaces of $\mathbb{R}^n$, and $\nu_m$ is the normalized Haar measure (that is, the unique rotationally invariant probability measure) on $G(n,m)$, while $\operatorname{vol}_m$ denotes $m$-dimensional volume. The limit case $q=0$ is understood in the logarithmic sense, namely by replacing the $q$-th power of the section volume with the logarithm of the section volume. For appropriate choices of $m$ and $q$, this family includes the
	dual quermassintegrals and the dual affine quermassintegrals as special
	cases.
	
	The centro-section measures introduced by Cai, Leng, Wu, and Xi
	\cite{2026-Cai-AdvMath} unify the dual curvature measures and their
	affine analogues. In the same work, they also established the
	corresponding variational formulas. For a convex body $K \subset \mathbb{R}^n$ which contains the origin in its interior, and a continuous function $f \in C(S^{n-1})$,
	\[
	\left.\frac{d}{dt}\right|_{t=0} \widetilde{\Psi}_{m,q}\big([K,f]_t\big) = q \int_{S^{n-1}} f(u) d\widetilde{A}_{m,q}(K,u).
	\] 
	The limit case $q=0$ has an analogous formulation based on the logarithmic definition of 
	$\widetilde{\Psi}_{m,0}.$ When $q=n$, this variational formula includes, as a special case, the
	variational formula for the affine dual quermassintegrals established in
	\cite{2025-Cai-AdvMath}.
	Here $[K,f]_t$ denotes the logarithmic family of Wulff shapes given by
	\[
	[K,f]_t = \big\{x \in \mathbb{R}^n : \log(x\cdot v)_+ \leq \log h_K(v) + t f(v),\ \text{for all } v \in S^{n-1}\big\},
	\]
	where $h_K$ is the support function of $K$, and $f\colon S^{n-1}\to\mathbb{R}$ is a continuous function. 
	
	The centro-section integrals studied in \cite{2026-Cai-AdvMath}
	depend on a single convex body. We introduce a mixed extension
	involving an independent reference body. In a mixed setting, however, it is natural to distinguish the body being varied from an independent reference body. Let $K,L\in\mathcal K_o^n$, $m\in\{1,\ldots,n-1\}$, and $q\in\mathbb R$. We introduce the $m$-th mixed centro-section integral of index $q$ by
	\begin{equation}
		\widetilde\Psi_{m,q}(K,L)
		=\int_{G(n,m)}
		\operatorname{vol}_m(K\cap\xi)^q\,
		\operatorname{vol}_m(L\cap\xi)^{n-q}
		\,d\nu_m(\xi). \label{1-2}
	\end{equation}
Although the functional $\widetilde{\Psi}_{m,q}(K,L)$ is well-defined for
$q\in\mathbb{R}$, the case $q=0$ is degenerate, since
$\widetilde{\Psi}_{m,0}(K,L)$ is independent of $K$. In contrast to the
logarithmic definition of the one-body functional in \eqref{1-1}, no
variation with respect to $K$ arises in this case. In the variational problem, we keep the reference body $L$ fixed
and vary $K$ through logarithmic Wulff shapes. For $q\neq0$, the mixed centro-section integral in \eqref{1-2}
	extends the $m$-th centro-section integrals of index
	$q$ studied in \cite{2026-Cai-AdvMath}; at $q=n$, it recovers the dual affine
	quermassintegral considered in \cite{2025-Cai-AdvMath}. Therefore, the framework includes the dual aﬃne quermassintegral and the $m$-th centro-section integrals of index $q$.  This also leads to a  new Minkowski-type problem : given a finite Borel measure on the unit sphere and a fixed reference body $L$, determine whether there exists a convex body $K$ such that its mixed centro-section measure equals the prescribed measure. 
	
	\textbf{Mixed centro-section measures.}  Let \(m\in\{1,\ldots,n-1\}\), \(q\in\mathbb R\), and \(K,L\in\mathcal K_o^n\). For a Borel set \(\eta\subset S^{n-1}\),  define
	\begin{equation}
		\widetilde A_{m,q}^{(L)}(K,\eta)
		=\int_{\boldsymbol{\alpha}_K^*(\eta)}\rho_K(u)^m
		\mathcal{R}_m^*\!\left(
		\operatorname{vol}_m(K\cap\cdot)^{q-1}
		\operatorname{vol}_m(L\cap\cdot)^{n-q}
		\right)(u)\,du.\label{1-3}
	\end{equation}
	where $\mathcal{R}_m^*$ is the dual Radon transform, and $\boldsymbol{\alpha}_K^*(\eta)$ denotes the reverse radial Gauss image of $K$, which is the set of all $u \in S^{n-1}$, such that the point $\rho_K(u)u$ has an outer unit normal to $\partial K$ that belongs to the set $\eta$. (Additional information is given in Section 2.)
	
	For $K,L\in  \mathcal{S}_o^n$, we define the mixed bi-dual $m$-intersection body $\mathbb{I}_{m,q}(K,L)$ of index $q$ to be the star body whose radial function is given by
	\begin{equation}
		\rho_{\mathbb{I}_{m,q}(K,L)}(u)^{n-m}
		= \mathcal{R}_m^*\Bigl(
		\mathrm{vol}_m(K\cap \cdot)^{q-1}
		\mathrm{vol}_m(L\cap \cdot)^{n-q}
		\Bigr)(u).\label{1-4}
	\end{equation}
	For a fixed convex body \( L \in \mathcal{K}_o^n \), let \( K \in \mathcal{K}_o^n \), and let \( K_t \) be a logarithmic family of Wulff shapes associated with a continuous function \( f \). We establish the following variational formula
	\begin{equation}\left.\frac{d}{dt}\right|_{t=0} \widetilde{\Psi}_{m,q}(K_t, L) = q \int_{S^{n-1}} f(v)\, d\widetilde{A}_{m,q}^{(L)}(K, v).\label{1-5}
	\end{equation}
	For fixed $L$, formula \eqref{1-5} expresses the first variation of the mixed centro-section integrals in terms of the mixed centro-section measure. If \(L=\omega_m^{-1/m}B^n\), the measure in \eqref{1-3} reduces to the centro-section measure in \cite{2026-Cai-AdvMath}, and \eqref{1-5} gives the corresponding variation formula. When $q=n$, the mixed centro-section integral is independent of $L$,
	and formula \eqref{1-5} yields the variational formula for dual affine
	quermassintegrals established in~\cite{2025-Cai-AdvMath}. We use this variational formula to prove the existence results stated below.

	\textbf{Minkowski problem for mixed centro-section measures.} Fix $m\in\{1,\dots,n-1\}$ and $q\in[0,n]$.
	For a given convex body $L\in\mathcal K_o^n$, the associated
	Minkowski problem is to characterize the finite Borel measures
	$\mu$ on $S^{n-1}$ for which there exists a convex body
	$K\in\mathcal K_o^n$ satisfying
	\[
	\widetilde{A}_{m,q}^{(L)}(K,\cdot)=\mu.
	\]
In this paper, we establish sufficient conditions for the solvability
of this problem.
	
	When $\mu$ has a density function $g:S^{n-1}\to\mathbb{R}$, the corresponding partial differential equation is a Monge–Ampère type equation on the sphere, involving Radon transforms $\mathcal{R}_m$, the dual Radon transform $\mathcal{R}_m^*$, and the spherical Hessian:
	\begin{align*}
		&\frac{1}{m^{n-1}}\mathcal{R}_m^*\Big(\mathcal{R}_m^{(q-1)} \circ \mathcal{D}^{(-m)}(h) \cdot \mathcal{R}_m^{(n-q)} \circ \mathcal{D}^{(-m)}(h_L)\Big)
		\left( \frac{\nabla h + h v}{|\nabla h + h v|} \right) \\
		&\quad \cdot h|\nabla h + h v|^{m-n} \det(\nabla^2 h + h I) = g, 
	\end{align*}
	where $h=h_K$ is the unknown support function of $K$, whereas $h_L$ is the known support function of the fixed reference body $L$. Here $\nabla h$ and $\nabla^2h$ denote the spherical gradient and spherical Hessian of \(h\), respectively. For \( h \in C^+(S^{n-1}) \), the
	\textit{duality transform} \( \mathcal{D} \) is defined as
	
	\[
	\mathcal{D}h(u) = \max_{v \in S^{n-1}} \frac{u \cdot v}{h(v)} \qquad \text{for } u \in S^{n-1}.
	\]
	Throughout, for an operator $\mathcal{A}$, when we write  $\mathcal{A}^{(k)}f(u)$, we shall mean $\mathcal{A}f(u) $ raised to the $k$-th power.
	
	\textbf{Existence for the mixed Minkowski problem.} As in the logarithmic Minkowski problem and its variants, the subspace concentration condition plays a critical role. A finite Borel measure $\mu$ on $S^{n-1}$ is said to satisfy the strict subspace concentration inequality if
	\begin{equation}
		\frac{\mu(S^{n-1}\cap\xi)}{|\mu|} < \frac{\dim\xi}{n}, \qquad \text{for each nonzero proper subspace } \xi\subset\mathbb{R}^n,\label{1-6}
	\end{equation}
where $|\mu|$ denotes the total mass of $\mu$. This condition was central to the solution of the even 
logarithmic Minkowski problem by Böröczky–Lutwak–Yang–Zhang \cite{2013-Boroczky-JAMS}.  Related necessary 
subspace concentration inequalities for cone-volume measures of origin-symmetric polytopes 
were proved independently by He–Leng–Li \cite{2006-He-AdvMath} and Henk–Schürmann–Wills  \cite{2005-Henk-Mathematika}.  In \cite{2025-Cai-AdvMath}, 
Cai, Leng, Wu, and Xi proved that \eqref{1-6} is suﬀicient for the aﬀine dual Minkowski problem;
in \cite{2026-Cai-AdvMath}, they established the corresponding existence result for centro-section measures. Our 
first main result extends these results to mixed centro-section measures.

	\begin{theorem}\label{thm-1-1} Let $m \in \{1,2,\dots,n-1\}$ be an integer. Suppose that $\mu$ is a finite even Borel measure on $S^{n-1}$ with positive mass and satisfies the strict subspace concentration inequality \eqref{1-6}. Then, for any fixed $L\in \mathcal{K}_o^n$, $0<q \leq n $, there exists an origin-symmetric convex body $K$ such that $$\widetilde{A}^{(L)}_{m,q}(K,\cdot) = \mu.$$
	\end{theorem}

Theorem~\ref{thm-1-1} recovers the corresponding existence results in \cite{2026-Cai-AdvMath} when \(L=\omega_m^{-1/m}B^n\), and in \cite{2025-Cai-AdvMath} when \(q=n\). For \(m=n-1\) and \(0<q<n\), the strict subspace concentration condition \eqref{1-6} can be replaced by the following weaker condition introduced in \cite{2026-Cai-AdvMath}:
	
	\begin{equation}
		\frac{\mu\big(S^{n-1} \cap \xi\big)}{|\mu|} < \frac{k}{n-1} - \frac{(k + q - n)_+}{q(n-1)},\quad \xi \in G(n,k),\ k = 1,\dots,n-1.
		\label{1-7}
	\end{equation}
	where $(a)_+ = \max\{a,0\}$ for $a \in \mathbb R$.
	
	\begin{theorem}\label{thm-1-2}
		Let $L\in \mathcal{K}_o^n$  be fixed. Suppose $q\in (0,n)$, and $\mu$ is a non-zero finite even Borel measure on $S^{n-1}$ that satisfies the refined subspace mass inequality \eqref{1-7}. Then there exists an origin-symmetric convex body $K \subset \mathbb R^n$ such that $\widetilde{A}^{(L)}_{n-1,q}(K,\cdot) = \mu$.
	\end{theorem}
	Theorems~\ref{thm-1-1} and~\ref{thm-1-2} hold for every fixed reference body \(L\in\mathcal K_o^n\), without any symmetry assumption on \(L\). The sufficient conditions on the prescribed measure are independent of the choice of \(L\). The proofs adapt the entropy and compactness estimates in \cite{2025-Cai-AdvMath,2026-Cai-AdvMath} to the weight determined by the central sections of \(L\). In this way, the existence results extend to mixed centro-section measures while retaining the subspace concentration hypotheses of the corresponding single-body problems.
	
	\textbf{ A sharp affine isoperimetric inequality.}
	In \cite{2025-Cai-AdvMath}, Cai, Leng, Wu, and Xi introduced
	the bi-dual $m$-intersection body operator $\mathbb{I}_m$
	and established a sharp affine isoperimetric inequality
	for its volume. In their subsequent work
	\cite{2026-Cai-AdvMath}, they introduced the bi-dual
	$m$-intersection bodies $\mathbb{I}_{m,q}K$ of index $q$
	associated with centro-section measures. However, the
	corresponding sharp volume inequality for $1<q<n$ was
	not established there.
	Motivated by these
	developments, we introduce the mixed bi-dual $m$-intersection body
	$\mathbb{I}_{m,q}(K,L)$ for a pair of star bodies and establish the
	following sharp affine isoperimetric inequality.
	
	\begin{theorem}\label{thm-1-3}
		\setlength{\parskip}{0pt}
		Let \(K,L\in\mathcal{S}_o^n\), let
		\(m\in\{1,\ldots,n-1\}\), and let \(1\leq q\leq n\). Then
		\vspace{-0.3em}
		\[
		V\bigl(\mathbb{I}_{m,q}(K,L)\bigr)
		\leq
		\left(
		\frac{m\omega_m^n}{n\omega_n^m}
		\right)^{\frac{n}{n-m}}
		V(K)^{\frac{m(q-1)}{n-m}}
		V(L)^{\frac{m(n-q)}{n-m}}.
		\]
		\vspace{-0.5em}
		
		If \(2\leq m\leq n-1\), the equality cases are as follows:
		\begin{enumerate}[
			label=(\alph*),
			leftmargin=2.8em,
			topsep=0.3em,
			itemsep=0.25em,
			parsep=0pt
			]
			\item If \(q=1\), equality holds if and only if \(L\) is an
			origin-centered ellipsoid, while \(K\) is arbitrary.
			
			\item If \(q=n\), equality holds if and only if \(K\) is an
			origin-centered ellipsoid, while \(L\) is arbitrary.
			
			\item If \(1<q<n\), equality holds if and only if \(K\) and \(L\)
			are origin-centered ellipsoids that are dilates of one another.
		\end{enumerate}
		\vspace{-0.3em}
		
		If \(m=1\), the same statements hold with
		``origin-centered ellipsoid'' replaced by
		``origin-symmetric star body.''
	\end{theorem}
	This result extends the affine isoperimetric inequality obtained in
	\cite{2025-Cai-AdvMath} to the mixed setting. Our proof combines the relation between the Radon transform and its dual with a dual mixed-volume representation of $\mathbb{I}_{m,q}(K,L)$. For
	$2\leq m\leq n-1$, the estimate follows from the generalized
	H\"older inequality and a sharp bound for the $n$-th moments of central
	section volumes. The case $m=1$ is treated separately by expressing the lengths of one-dimensional central sections as sums of radial functions in opposite directions and applying H\"older's inequality together with a sharp elementary convexity inequality. 
	
	The paper is organized as follows. Section 2 collects the necessary background on convex bodies, radial Gauss images, and the extended Radon transform. Section 3 introduces the mixed centro-section measures and establishes their variational formulas. Section 4 reformulates the Minkowski problem as a maximization problem for a suitable functional involving entropy. Section 5 proves the existence results under the subspace concentration condition and a refined subspace mass inequality, thereby establishing Theorems~\ref{thm-1-1} and~\ref{thm-1-2}. Section 6 studies the affine transformation properties of the mixed centro-section integrals and the mixed bi-dual $m$-intersection body operator $\mathbb{I}_{m,q}(K,L)$, and proves the affine isoperimetric inequality stated in Theorem~\ref{thm-1-3}.
	
	\section{Preliminaries}
	
	\qquad We recall some basic notions from the theory of convex bodies and integral
	geometry that will be used below. Details can be found in \cite{2014-Schneider-CUP} and \cite{2006-Gardner-CUP}.
	
	We work in the $n$-dimensional Euclidean space $\mathbb{R}^n$, equipped with the canonical inner
	product $x\cdot y$ and the norm $|x|=\sqrt{x\cdot x}$. The $m$-dimensional Hausdorff measure is denoted
	by $\mathcal{H}^m$. When a convex body $K$ lies in an $m$-dimensional affine subspace, its m-dimensional volume is denoted by $\operatorname{vol}_m(K)$  rather than $\mathcal{H}^m(K)$. The $m$-dimensional volume of the $m$-dimensional unit ball is denoted by $\omega_m$.
	For a finite measure $\mu$ on $S^{n-1}$, we write $|\mu|$ for its total mass, that is $|\mu|=\mu(S^{n-1})$.
	
	\subsection*{2.1.\ Convex bodies and star bodies}
	\quad   The set of convex bodies in $\mathbb{R}^n$ is denoted by $\mathcal{K}^n$, the set of convex bodies containing the origin in their interior is denoted by $\mathcal{K}_o^n$, and $\mathcal{K}_e^n\subset\mathcal{K}_o^n$
	will be used to denote the set of origin-symmetric convex bodies.
	
	If $K \subset \mathbb{R}^n$ is a nonempty compact convex set, then its support function is defined by
	$$
	h_K(x) = \max\{ x \cdot y : y \in K \}.
$$
	The support function is convex and positively homogeneous of degree one. It is immediate from the definition that, for any $\phi \in \mathrm{SL}(n)$ and any $x \in \mathbb{R}^n$, the support function of the image $\phi K$ satisfies
	\[
	h(\phi K, x) = h(K, \phi^t x),
	\]
	where $\phi^t$ is the transpose of $\phi$. Moreover, the support function is continuous on $S^{n-1}$ and uniquely determines the convex body $K \subset \mathbb{R}^n$ by $x \cdot v \leq h(v)$ for all $v \in S^{n-1}$.
	
	The collection of all nonempty compact convex sets is viewed as a metric space with the Hausdorff metric,
	and the distance between $K,L\in\mathcal{K}^n$ is the max-norm of the absolute value of the difference of their
	support functions, $\|h_K-h_L\|_\infty$.
	
	For a set $K\subset\mathbb{R}^n$ star-shaped with respect to the origin, its radial function $\rho_K:S^{n-1}\to[0,\infty),$ is
	defined by
	\[
	\rho_K(u)=\max\{t\geq 0:\ tu\in K\}.
	\]
	The radial function can be extended to
	\(\mathbb R^n\setminus\{0\}\) as a function homogeneous of degree \(-1\). A star-shaped set $K$ is said to be a \emph{star body} if $\rho_K$ is continuous and positive. The collection of star bodies is denoted by $\mathcal{S}_o^n$. Also note that the boundary of \( K \), written as \( \partial K \), is defined by
	$$\partial K = \left\{ \rho_K(u)u : u \in S^{n-1} \right\}.$$
	
	Let $q\ne 0$. Given $K,L\in\mathcal{S}^n_o$, the \emph{dual mixed volume} is defined by
	\[
	\widetilde{V}_q(K,L)=\frac{1}{n}\int_{S^{n-1}} \rho_K(u)^{n-q}\rho_L(u)^q\,du.
	\]
	Dual mixed volumes are affine invariant under simultaneous $SL(n)$-transformations of
	the bodies; for each $\varphi\in SL(n)$,
	\[
	\widetilde{V}_q(\varphi K,\varphi L)=\widetilde{V}_q(K,L).
	\]
	For convenience, define the $m$-th dual volume
	\[
	\widetilde{V}_m(K) = \widetilde{W}_{n-m}(K) = \widetilde{V}_m(B^n, K) = \frac{1}{n} \int_{S^{n-1}} \rho_K(u)^m \, du,
	\]
	where $B^n$ denotes the unit ball in $\mathbb{R}^n$. Applying H\"{o}lder's inequality, one obtains the dual mixed volume inequality,
	\begin{equation}
		\widetilde{V}_m(K, L) \leq V(K)^{(n-m)/n} V(L)^{m/n}, \label{2-1}
	\end{equation}
	with equality if and only if $K$ and $L$ are dilates.
	
	\subsection*{2.2.\ The radial Gauss image and its reverse} \qquad For $K\in\mathcal{K}_o^n$ and $\omega\subset S^{n-1}$, the \emph{radial Gauss image} $\alpha_K(\omega)$ is the set of all outer unit
	normals of $K$ at the boundary point $\rho_K(u)u$, where $u\in\omega$. That is,
	\begin{equation*}
		\boldsymbol{\alpha}_K(\omega)=\bigcup_{u\in\omega}\left\{v\in S^{n-1}:\ \rho_K(u)u\cdot v=h_K(v)\right\}.
	\end{equation*}
	\qquad The \emph{reverse radial Gauss image} $\boldsymbol{\alpha}_K^*(\eta)$ is the set of all radial directions $u\in S^{n-1}$,
	such that $\rho_K(u)u$ has at least one element in $\eta$ as outer unit normal, that is
	\begin{equation*}
		\boldsymbol{\alpha}_K^*(\eta)=\bigcup_{v\in\eta}\left\{u\in S^{n-1}:\ \rho_K(u)u\cdot v=h_K(v)\right\}.
	\end{equation*}
	
	Let $\omega_K\subset S^{n-1}$ denote the set of $u\in S^{n-1}$ such that $\boldsymbol{\alpha}_K(\{u\})$ is not a singleton, that
	is $\rho_K(u)u$ has more than one outer unit normal. The \emph{radial Gauss map}
	\[
	\alpha_K:S^{n-1}\setminus\omega_K\longrightarrow S^{n-1},
	\]
	is defined by $\boldsymbol{\alpha}_K(\{u\})=\{\alpha_K(u)\}$. It can be shown that $\mathcal{H}^{n-1}(\omega_K)=0$ (see \cite{2014-Schneider-CUP}).
	
	If $\{\eta_j\}_{j=1}^\infty$ is a sequence of pairwise disjoint sets in $S^{n-1}$, then
	\begin{equation*}
		\{\boldsymbol{\alpha}_K^*(\eta_j)\setminus\omega_K\}_{j=1}^\infty \ \text{is pairwise disjoint as well.}
	\end{equation*}
	
	Let $\eta_{K}\subset S^{n-1}$ denote the set of all $v\in S^{n-1}$ for which $\boldsymbol{\alpha}_K^*(\{v\})$ is not a singleton,
that is, there is more than one radial direction $u\in S^{n-1}$ such that $v$ is the unit outer
	normal of $\rho_K(u)u$. The \emph{reverse radial Gauss map},
	\[
	{\alpha}_K^*:S^{n-1}\setminus\eta_K\longrightarrow S^{n-1},
	\]
	is defined, for $u\in S^{n-1}\setminus\eta_K$, by $\boldsymbol{\alpha}_K^*(\{u\})=\{\alpha_K^*(u)\}$. The set $\eta_K$ is also of $\mathcal{H}^{n-1}$-measure
	$0$ (see Theorem 2.2.11 in Schneider \cite{2014-Schneider-CUP}).
	
	\subsection*{2.3.\ The extended Radon transform and its dual}
	\qquad Let $m\in\{1,\ldots,n-1\}$. The $m$-dimensional Radon transform, $\mathcal{R}_m:C(S^{n-1})\to C(G(n,m))$,
	is the linear operator defined by
	\[
	\mathcal{R}_m f(\xi)=\int_{S^{n-1}\cap \xi} f(w)\,dw,
	\]
	where the integration is with respect to $(m-1)$-dimensional spherical Lebesgue measure.
	The $m$-dimensional dual Radon transform, $\mathcal{R}_m^*:C(G(n,m))\to C(S^{n-1})$ is
	defined by
	\begin{equation}\label{2-2}
		\mathcal{R}_m^*F(u)=\frac{m\omega_m}{n\omega_n}\int_{G_u(n-1,m-1)}
		F(\operatorname{span}\{u,\zeta\})\,d\nu_{m-1}(\zeta),
	\end{equation}
	where $\nu_{m-1}$ denotes Haar measure on $G_u(n-1,m-1)$, the Grassmannian of $(m-1)$-
	dimensional subspaces of $u^\perp$.
	
	The Radon transform $\mathcal{R}_m$ and its dual are related by (see, e.g. \cite{2005-Koldobsky-AMS}),
	\[
	\int_{G(n,m)} F(\xi)\mathcal{R}_m f(\xi)\,d\nu_m(\xi)
	=
	\int_{S^{n-1}} f(u)\mathcal{R}_m^*F(u)du \quad \text{for}
	\; f\in C(S^{n-1}),\;\ F\in C(G(n,m))\quad
	\]
	which we write as,
	\begin{equation*}
		(F,\mathcal{R}_m f)=(\mathcal{R}_m^*F,f) \quad \text{for} \; f\in C(S^{n-1}),\quad F\in C(G(n,m)).
	\end{equation*}
	A direct consequence is that
	\begin{equation*}
		\int_{G(n,m)}\left(\int_{S^{n-1}\cap \xi} f(u)\,du\right)d\nu_m(\xi)
		=\frac{m\omega_m}{n\omega_n}\int_{S^{n-1}} f(u)\,du,
	\end{equation*}
	for $f\in C(S^{n-1})$, and where the inner integral on the left is with respect to $(m-1)$-
	dimensional spherical Lebesgue measure and the integration on the right is with respect
	to $(n-1)$-dimensional spherical Lebesgue measure.
	
	When $m=n-1$, since $\xi\in G(n,n-1)$ can be represented by $u^\perp$ for $u\in S^{n-1}$, the
	transform $\mathcal{R}_{n-1}$ can be identified with the \emph{spherical Radon transform}
	$\mathcal{R}:C(S^{n-1})\to C(S^{n-1})$, which is defined by
	\begin{equation*}
		\mathcal{R}f(u)=\int_{S^{n-1}\cap u^\perp} f(w)\,dw,
	\end{equation*}
	where the integration is with respect to $(n-2)$-dimensional spherical Lebesgue measure.
	The operator $\mathcal{R}$ is self-adjoint in that for all $f,g\in C(S^{n-1})$, we have
	$(\mathcal{R}f,g)=(f,\mathcal{R}g)$, or equivalently,
	\[
	\int_{S^{n-1}} g(u)\mathcal{R}f(u)\,du=\int_{S^{n-1}} f(u)\mathcal{R}g(u)\,du,
	\]
	where the integration is with respect to $(n-1)$-dimensional spherical Lebesgue measure.
	
	The classical \emph{intersection body} $\mathrm{I}K$ of a convex body $K\in\mathcal{K}_o^n$ introduced by Lutwak \cite{1988-Lutwak-AdvMath}, is defined via the Radon transform by
	$$\rho_{\mathrm{I}K}(u)=\operatorname{vol}_{n-1}(K\cap u^\perp)=\mathcal{R}\left(\rho_K^{\,n-1}/(n-1)\right)(u).
	$$	
	As shown in \cite[Theorem 3.2, Theorem 3.3]{2025-Cai-AdvMath}, for each $p\in[1,+\infty]$, the $m$-dimensional Radon transform and its dual extend to the relevant $L^{p}$ spaces. The extended Radon transform $\mathcal{R}_m$ is a bounded linear operator from
	$L^p(S^{n-1})$ to $L^p(G(n,m))$, which satisfies
	\begin{equation}\label{2-3}
		\|\mathcal{R}_m f\|_p\le m\omega_m\|f\|_p, \qquad f\in L^p(S^{n-1}),
	\end{equation}
	where $\|\cdot\|_p$ denotes the $L^p$ norm. The dual transform $\mathcal{R}_m^*$ is also a bounded linear operator from $L^p(G(n,m))$ to
	$L^p(S^{n-1})$, which satisfies
	\begin{equation}\label{2-4}
		\|\mathcal{R}_m^*F\|_p\le m\omega_m\|F\|_p, \qquad F\in L^p(G(n,m)) .
	\end{equation}
	
	\begin{theorem}[{\cite[Theorem 3.4]{2025-Cai-AdvMath}}]
		\label{thm-2-1}
		Let $p,p'\in[1,+\infty]$ satisfy
		$
		\frac{1}{p}+\frac{1}{p'}=1.
		$
		For $f\in L^p(S^{n-1})$ and $F\in L^{p'}(G(n,m))$, the extended Radon transforms satisfy the duality relation
		\begin{equation}\label{2-5}
			\int_{G(n,m)}
			\mathcal{R}_m f(\xi)\,F(\xi)\,d\nu_m(\xi)
			=
			\int_{S^{n-1}}
			f(u)\,\mathcal{R}_m^{*}F(u)\,du.
		\end{equation}
		where the integral on the right-hand side is taken with respect to the $(n-1)$-dimensional spherical Lebesgue measure.
	\end{theorem}
	
	\section{General mixed centro-section measures and the variation formula}
	
	In this section, we introduce mixed centro-section measures and
	establish the corresponding variational formula. This new family
	contains several known geometric measures as special cases.
	
	\subsection*{3.1\ Mixed centro-section measures}
	\noindent\textbf{Definition 3.1.}
	Suppose $q\in\mathbb{R}$ and $m\in\{1,2,\ldots,n-1\}$.
	Let $K,L\in\mathcal{K}_o^n$, and $\eta\subset S^{n-1}$ be a Borel set. We define
	\begin{equation}\label{3-1}
		\widetilde{A}^{(L)}_{m,q}(K,\eta)
		=
		\int_{\boldsymbol{\alpha}_K^*(\eta)} \rho_K(u)^m\,
		\mathcal{R}_m^{*}\bigl(\mathrm{vol}_m(K\cap \cdot)^{\,q-1}\mathrm{vol}_m(L\cap \cdot)^{n-q}\bigr)(u)du,
	\end{equation}
	where the inner integral on the right is defined by \eqref{2-2}, while the outer integration is
	with respect to $(n-1)$-dimensional spherical Lebesgue measure. Utilizing the notion of a mixed bi-dual $m$-intersection body \eqref{1-4}, $\widetilde{A}^{(L)}_{m,q}(K,\eta)$ can be rewritten as
	\begin{equation*}
		\widetilde{A}_{m,q}^{(L)}(K,\eta) = \int_{\boldsymbol{\alpha}_K^*(\eta)}\rho_K(u)^m\rho_{\mathbb{I}_{m,q}(K,L)}(u)^{n-m}\,du,
	\end{equation*}
	where the integration is with respect to spherical $(n-1)$-dimensional Lebesgue measure.
	
	We next describe the relations between the mixed centro-section measure
	introduced above and several known geometric measures.
	
	\setcounter{theorem}{1}
	\noindent\textbf{Remark 3.2.} 
		\textnormal{ \textbf{(1) Relation with generalized dual curvature measures.}}
		The generalized dual curvature measures introduced by Lutwak, Yang,
		and Zhang \cite{2018-Lutwak-Advmath} are defined, for
		$K\in\mathcal K_o^n,Q\in\mathcal S_o^n$ and a Borel set
		$\eta\subset S^{n-1}$, by
		\[
		\widetilde C_m(K,Q,\eta)
		=
		\int_{\boldsymbol{\alpha}_K^*(\eta)}
		\rho_K(u)^m
		\rho_Q(u)^{n-m}\,du,
		\]
		where we use \(n\) times the normalization in \cite{2018-Lutwak-Advmath} .By the definition \eqref{1-4} of the mixed bi-dual \(m\)-intersection body
		\(\mathbb I_{m,q}(K,L)\), we have
		\[
		\widetilde A_{m,q}^{(L)}(K,\eta)
		=
		\widetilde C_m
		\bigl(
		K,\mathbb I_{m,q}(K,L),\eta
		\bigr).
		\]
		Hence, the mixed centro-section measure can be regarded as a generalized
		dual curvature measure whose second argument is determined by the mixed
		centro-section integrals of the pair \((K,L)\).
		
		\textnormal{ \textbf{(2) Reduction to centro-section and affine dual curvature measures.}}
		If \(L=B^n\), then
		$
		\operatorname{vol}_m(B^n\cap\xi)^{n-q}
		=
		\omega_m^{n-q},
		$
		for every \(\xi\in G(n,m)\). Hence,
		\[
		\begin{aligned}
			\widetilde A_{m,q}^{(B^n)}(K,\eta)
			&=
			\omega_m^{n-q}
			\int_{\boldsymbol{\alpha}_K^*(\eta)}
			\rho_K(u)^m
			\mathcal R_m^*
			\left(
			\operatorname{vol}_m(K\cap\cdot)^{q-1}
			\right)(u)\,du
			\\
			&=
			\omega_m^{n-q}
			\widetilde A_{m,q}(K,\eta).
		\end{aligned}
		\]
		Moreover, if either \(q=n\) or \(K=L\), then
		$
		\operatorname{vol}_m(K\cap\xi)^{q-1}
		\operatorname{vol}_m(L\cap\xi)^{n-q}
		=
		\operatorname{vol}_m(K\cap\xi)^{n-1}.
		$
		Consequently,
		$
		\widetilde A_{m,q}^{(L)}(K,\eta)
		=
		\widetilde C_m^a(K,\eta)
		$
		when \(q=n\), while
		$
		\widetilde A_{m,q}^{(K)}(K,\eta)
		=
		\widetilde C_m^a(K,\eta)
		$
		when \(K=L\). Thus, the mixed centro-section measures contain both the centro-section measures and the affine dual curvature measures as special cases.
		
		\textnormal{ \textbf{(3) The hyperplane section case.}}
		If \(m=n-1\), then the mixed bi-dual $(n-1)$-intersection body can be represented in terms of the classical intersection bodies $	\rho_{\mathrm I K} $ and $\rho_{\mathrm I L}.$ In particular,
		$
		\rho_{\mathrm I K}(u)
		=
		\operatorname{vol}_{n-1}(K\cap u^\perp),
		\;
		\rho_{\mathrm I L}(u)
		=
		\operatorname{vol}_{n-1}(L\cap u^\perp).
		$
		In this case, \(\mathcal R_{n-1}\) is naturally identified with the
		spherical Radon transform \(\mathcal R\). Moreover, under the normalization
		used for the dual Radon transform,
		$
		\mathcal R_{n-1}^{*}
		=
		\frac1{n\omega_n}\mathcal R .
		$
		Therefore,
		\[
		\widetilde A_{n-1,q}^{(L)}(K,\eta)
		=
		\frac1{n\omega_n}
		\int_{\boldsymbol{\alpha}_K^*(\eta)}
		\rho_K(u)^{n-1}
		\mathcal R
		\left(
		\rho_{\mathrm I K}^{q-1}
		\rho_{\mathrm I L}^{n-q}
		\right)(u)\,du .
		\]
		In particular, when \(q=n\), the dependence on \(L\) disappears. Since
		\[
		\mathcal R
		\left(
		\rho_{\mathrm I K}^{n-1}
		\right)(u)
		=
		(n-1)\rho_{\mathrm I^2K}(u),
		\]
		we obtain
		\[
		\widetilde A_{n-1,n}^{(L)}(K,\eta)
		=
		\frac{n-1}{n\omega_n}
		\int_{\boldsymbol{\alpha}_K^*(\eta)}
		\rho_K(u)^{n-1}
		\rho_{\mathrm I^2K}(u)\,du .
		\]
		Thus, in the hyperplane section case with \(q=n\), the mixed
		centro-section measure is independent of \(L\) and coincides with the
		affine dual curvature measure up to a normalization constant.
		
		\textnormal{\textbf{(4) Finiteness and the push-forward representation.}}
		Taking \(\eta=S^{n-1}\), it follows from \eqref{2-5}
		that
		\[
		\begin{aligned}
			\widetilde A_{m,q}^{(L)}
			(K,S^{n-1})
			&=
			m
			\int_{G(n,m)}
			\operatorname{vol}_m(K\cap\xi)^q
			\operatorname{vol}_m(L\cap\xi)^{n-q}
			\,d\nu_m(\xi)=
			m\widetilde\Psi_{m,q}(K,L)
			<\infty .
		\end{aligned}
		\]
		Therefore,
		$
		\widetilde A_{m,q}^{(L)}(K,\cdot)
		$
		is a finite Borel measure on $S^{n-1}.$ Moreover, it is even whenever
		\(K\in\mathcal K_e^n.\) Finally, since $\mathcal H^{n-1}(\omega_K)=0$, we may regard $\alpha_K$ as a Borel measurable map on $S^{n-1}$ by assigning arbitrary values on $\omega_K$. The choice of these values does not affect the integrals below. The measure $\widetilde A_{m,q}^{(L)}(K,\cdot)$ is the push-forward under $\alpha_K$ of the measure that is absolutely continuous with respect to $\mathcal H^{n-1}$ and has density
		\[
		g_{K,L}(u)
		=
		\rho_K(u)^m
		\mathcal R_m^*
		\left(
		\operatorname{vol}_m(K\cap\cdot)^{q-1}
		\operatorname{vol}_m(L\cap\cdot)^{n-q}
		\right)(u).
		\]
		Equivalently,
		\[
		\widetilde A_{m,q}^{(L)}(K,\cdot)
		=
		(\alpha_K)_{\#}
		\bigl(g_{K,L}\mathcal H^{n-1}\bigr).
		\]
		Consequently, by the integration formula for push-forward measures, for every bounded Borel function $f:S^{n-1}\to\mathbb R$,
		\[
		\int_{S^{n-1}}
		f(v)\,d\widetilde A_{m,q}^{(L)}(K,v)
		=
		\int_{S^{n-1}}
		f\bigl(\alpha_K(u)\bigr)g_{K,L}(u)\,
		d\mathcal H^{n-1}(u).
		\]
		That is,
		\begin{equation}\label{3-2}
			\begin{aligned}
				\int_{S^{n-1}}
				f(v)\,d\widetilde A_{m,q}^{(L)}(K,v)
				&=
				\int_{S^{n-1}}
				f\bigl(\alpha_K(u)\bigr)\rho_K(u)^m
				\\
				&\quad\cdot
				\mathcal R_m^*
				\left(
				\operatorname{vol}_m(K\cap\cdot)^{q-1}
				\operatorname{vol}_m(L\cap\cdot)^{n-q}
				\right)(u)\,
				d\mathcal H^{n-1}(u).
			\end{aligned}
			\tag{3.2}
		\end{equation}

	The mixed centro-section measures arise as the first variations of the mixed centro-section integrals.
	Define
	\begin{equation}\label{3-3}
		\widetilde{\Psi}_{m,q}(K,L)
		=
		\int_{G(n,m)}
		\mathrm{vol}_m(K\cap \xi)^{q}
		\mathrm{vol}_m(L\cap \xi)^{n-q}
		\,d\nu_m(\xi).
		\tag{3.3}
	\end{equation}
	where \(m\in\{1,2,\ldots,n-1\}\) is an integer. 
	
	\noindent\textbf{Remark 3.3.} 
	For $K,L\in\mathcal{K}_o^n$, the integral in (3.3) is
	finite for every $q\in\mathbb{R}$. If the origin lies
	on the boundary of $K$ and $q<0$, it may be infinite.
	For fixed $L\in\mathcal{K}_o^n$ and
	$K_1,K_2\in\mathcal{K}_o^n$ with $K_1\subsetneq K_2$,
	it is strictly increasing in $K$ when $q>0$ and
	strictly decreasing when $q<0$:
	\[
	\begin{aligned}
		\widetilde{\Psi}_{m,q}(K_1,L)
		&<
		\widetilde{\Psi}_{m,q}(K_2,L)\qquad(q>0),\\
		\widetilde{\Psi}_{m,q}(K_1,L)
		&>
		\widetilde{\Psi}_{m,q}(K_2,L)\qquad(q<0).
	\end{aligned}
	\]
	For $q=0$, it is independent of $K$.
	
	\textbf{Logarithmic family of Wulff shapes.}
	Let $\Omega\subset S^{n-1}$ be a closed set not contained in any closed hemisphere of $S^{n-1}$.
	Let $h_0:\Omega\to(0,\infty)$ be continuous and let $f:\Omega\to\mathbb R$ be continuous.
	For each $t\in\mathbb R$, define $h_t:\Omega\to(0,\infty)$ by
	\[
	\log h_t(v)=\log h_0(v)+t f(v),\qquad v\in\Omega .
	\]
	The logarithmic family of Wulff shapes associated with $(h_0,f)$ is defined by
	\[
	K_t=\{x\in\mathbb R^n:x\cdot v\le h_t(v)\ \text{for all }v\in\Omega\},\qquad t\in\mathbb R .
	\]
	We write $K=K_0$. In particular, $K_t\in\mathcal K_o^n$ for all $t$.
	Throughout this section, let $L\in\mathcal K_o^n$ be fixed.
	
	\medskip
	\begingroup
	\renewcommand{\thetheorem}{3.4}
	\begin{theorem}\label{thm-3-4}
		Suppose that $0< q \leq n $ and \(m\in\{1,2,\ldots,n-1\}\).
		Let \(K_t\) be the logarithmic family of Wulff shapes associated with \((h_0,f)\),
		and write \(K=K_0\). Let \(L\in\mathcal K_o^n\) be fixed. 
		Then,
		\begin{equation}\label{3-4}
			\lim_{t\to 0}\frac{\widetilde{\Psi}_{m,q}(K_t,L)-\widetilde{\Psi}_{m,q}(K,L)}{t}
			=
			q\int_{S^{n-1}} f(v)\, d\widetilde{A}^{(L)}_{m,q}(K,v),
			\tag{3.4}
		\end{equation}
	\end{theorem}
	\noindent\textbf{Remark 3.5.} If $q=0$, then
	\[
	\widetilde{\Psi}_{m,0}(K,L)
	=
	\int_{G(n,m)}
	\operatorname{vol}_m(L\cap\xi)^n\,d\nu_m(\xi)
	=
	\widetilde{\Psi}_{m,n}(L),
	\]
	which is independent of $K$. Consequently, for any logarithmic family
	$\{K_t\}$ with $K_0=K$,
	\[
	\left.
	\frac{d}{dt}
	\right|_{t=0}
	\widetilde{\Psi}_{m,0}(K_t,L)=0.
	\]
	Thus, when $q=0$, the first variation with respect to $K$ is identically zero. Consequently, the variational formula does not yield a nontrivial Minkowski problem in this case.

	\endgroup
	
	\subsection*{3.2\ Variation formulas for mixed centro-section measures}
	
	We are now in a position to prove Theorem~\ref{thm-3-4}, and we will use the following lemma,
	which is a direct corollary of Lemma 4.2 in \cite{2016-Huang-ActaMath}.
	
	\noindent\textbf{Lemma 3.6.}
	\begingroup
	\makeatletter
	\phantomsection
	\def\@currentlabel{3.6}
	\label{lem-3-6}
	\makeatother
	\endgroup
	\textit{\cite[Lemma 3.8]{2025-Cai-AdvMath}}
	Let $K_t$ be the logarithmic family of Wulff shapes defined by
	\begin{equation}	\label{3-5}
		K_t=\{x\in\mathbb{R}^n:\ x\cdot v\le h_t(v)\ \text{for all }v\in\Omega\},
		\qquad t\in\mathbb{R}.
		\tag{3.5}
	\end{equation}
	Then
	\[
	\mathcal{R}_m\!\left(
	\frac{\rho_{K_t}^{\,m}-\rho_K^{\,m}}{t}
	\right)
	\to
	m\mathcal{R}_m\bigl(\rho_K^{\,m}\cdot f\circ\alpha_K\bigr),
	\qquad \text{in }L^p(G(n,m)),
	\]
	for every $p\in[1,\infty)$. Moreover, there exist $\delta_0>0$ and $M_0>0$ so that
	\begin{equation}	\label{3-6}
		\bigl|\rho_{K_t}(u)^m-\rho_K(u)^m\bigr|
		\le M_0|t|,
		\qquad |t|<\delta_0,\ u\in\Omega.
		\tag{3.6}
	\end{equation}
	
	\medskip
	\begin{proof}[Proof of Theorem~\ref{thm-3-4}]
		Using Lemma~\ref{lem-3-6}, as $t\to 0$,
		\[
		\mathcal{R}_m\left(\frac{\rho_{K_t}^m - \rho_K^m}{t}\right)\mathrm{vol}_m(L\cap\xi)^{n-q} \to m\mathcal{R}_m\big(\rho_K^m\cdot f\circ\alpha_K\big)\mathrm{vol}_m(L\cap\xi)^{n-q},\quad \text{in }L^p(G(n,m)),
		\]
		for each $p\in[1,\infty)$. Furthermore, from \eqref{3-6} and the boundedness estimate \eqref{2-3}, we have
		\begin{equation}	\label{3-7}
			\big|\mathcal{R}_m\big(\rho_{K_t}^m - \rho_K^m\big)\big|\mathrm{vol}_m(L\cap\xi)^{n-q} \leq m\omega_m M_0|t|\mathrm{vol}_m(L\cap\xi)^{n-q},
			\tag{3.7}
		\end{equation}
		for all $\xi\in G(n,m)$ and $t\in(-\delta_0,\delta_0)$. Let $(t_i)$ be an arbitrary sequence in $\mathbb{R}\setminus\{0\}$ that converges to $0$ as $i\to\infty$. Recall that
		\[
		\mathrm{vol}_m(K_{t_i}\cap\xi) = \frac{1}{m}\int_{S^{n-1}\cap\xi}\rho_{K_{t_i}}(u)^m du = \frac{1}{m}\cdot\mathcal{R}_m(\rho_{K_{t_i}}^m)(\xi).
		\]
		Using Lemma~\ref{lem-3-6}, we conclude that there exists a subsequence $(t_{i_j})$ such that
		\begin{equation}	\label{3-8}
			\frac{\mathrm{vol}_m(K_{t_{i_j}}\cap\xi) - \mathrm{vol}_m(K\cap\xi)}{t_{i_j}} \to \mathcal{R}_m\big(\rho_K^m\cdot f\circ\alpha_K\big)(\xi),\quad \text{as }j\to\infty,\text{ a.e. }\xi\in G(n,m).
			\tag{3.8}
		\end{equation}
		Suppose $q\neq 0$. The differential mean value theorem tells us that for $0<a<b$,
		\[
		|a^q - b^q| = q\theta^{q-1}|a - b|,\quad \text{for some }\theta\in(a,b).
		\]
		Since $\rho_{K_t}(u)$ is uniformly bounded from both above and below for $t\in(-\delta_0,\delta_0)$, there exists a uniform constant $c_0>0$, such that
		\[
		q\cdot\big(\mathrm{vol}_m(K_t\cap\xi)\big)^{q-1} \leq c_0.
		\]
		As a result, by these and the previous inequality  \eqref{3-6},
		\[
		\begin{aligned}
			&\left|
			\frac{\big(\mathrm{vol}_m(K_{t}\cap\xi)^{q} - \mathrm{vol}_m(K\cap\xi)^q\big)\mathrm{vol}_m(L\cap\xi)^{n-q}}{t}
			\right| \\
			&\leq c_0\cdot\left|\frac{\mathrm{vol}_m(K_t\cap\xi) - \mathrm{vol}_m(K\cap\xi)}{t}\right|\mathrm{vol}_m(L\cap\xi)^{n-q} \\
			&\leq c_0\cdot\omega_m\cdot M_0\mathrm{vol}_m(L\cap\xi)^{n-q}.
		\end{aligned}
		\]
		Therefore, by \eqref{3-8} and the dominated convergence theorem, and then by the duality of extended Radon transform \eqref{2-5} and \eqref{3-2},  we obtain
		\[
		\begin{aligned}
			&\lim_{j\to\infty} \frac{\widetilde{\Psi}_{m,q}(K_{t_{i_j}},L) - \widetilde{\Psi}_{m,q}(K,L)}{t_{i_j}} \\
			&\quad= \lim_{j\to\infty}\int_{G(n,m)} \frac{\big(\mathrm{vol}_m(K_{t_{i_j}}\cap\xi)^q - \mathrm{vol}_m(K\cap\xi)^q\big)\mathrm{vol}_m(L\cap\xi)^{n-q}}{t_{i_j}}d\nu_m(\xi)\\
			&\quad= q\int_{G(n,m)} \mathrm{vol}_m(K\cap\xi)^{q-1}\mathrm{vol}_m(L\cap\xi)^{n-q}\cdot \mathcal{R}_m\big(\rho_K^m\cdot f\circ\alpha_K\big)(\xi)\,d\nu_m(\xi)\\
			&\quad= q\int_{S^{n-1}} \rho_K(u)^m f\big(\alpha_K(u)\big) \mathcal{R}_m^*\big(\mathrm{vol}_m(K\cap \cdot)^{q-1}\mathrm{vol}_m(L\cap\cdot)^{n-q}\big)(u)\,du\\
			&\quad= q\int_{S^{n-1}} f(v)\,d\widetilde{A}_{m,q}^{(L)}(K,v).
		\end{aligned}
		\]
		Since $(t_i)$ is an arbitrary sequence converging to $0$, we obtain that
		\[
		\lim_{t\to 0}\frac{\widetilde{\Psi}_{m,q}(K_t,L) - \widetilde{\Psi}_{m,q}(K,L)}{t}
		= q\int_{S^{n-1}} f(v)\,d\widetilde{A}_{m,q}^{(L)}(K,v).
		\]
		This proves \eqref{3-4}.\\
		
	\end{proof}
	
	\section{Solution to a maximization problem}	
	\quad Let $\mu$ be a finite even Borel measure on $S^{n-1}$ with $|\mu| = \mu(S^{n-1}) > 0.$
	For $K \in \mathcal{K}_e^n$, define the entropy functional
	\[
	E_\mu(K) = -\frac{1}{|\mu|}\int_{S^{n-1}} \log h_K(v)\,d\mu(v).
	\]
	Fix $L \in \mathcal{K}_o^n$ and $q \in (0,n]$. For $K \in \mathcal{K}_e^n$, set
	\[
	J_{m,q}^{(L)}(K) = \frac{1}{mq}\log \widetilde{\Psi}_{m,q}(K,L) + E_\mu(K).
	\]
	This functional is homogeneous of degree zero with respect to the body \(K\).
	More precisely, for every \(c>0\),
	$
	J_{m,q}^{(L)}(cK)=J_{m,q}^{(L)}(K).
	$
	In particular, when \(q=n\), the functional is independent of \(L\), and
	\[
	J_{m,n}^{(L)}(B^n)
	=
	\frac{1}{m}\log\omega_m .
	\]
	\noindent\textbf{Maximization Problem.}
	Given a finite even Borel measure 
	\(\mu\) on \(S^{n-1}\) and a fixed convex body \(L\) containing the origin 
	in its interior, we consider the following maximization problem: determine 
	whether there exists \(K_0\in\mathcal K_e^n\) such that
	\begin{equation}\label{4-1}
		J_{m,q}^{(L)}(K_0)
		=
		\sup\left\{
		J_{m,q}^{(L)}(K):K\in\mathcal K_e^n
		\right\}.
	\end{equation}
	
	We shall show that the solution of the Maximization Problem is the solution of the
	Minkowski problem for the mixed centro-section measure $\widetilde{A}^{(L)}_{m,q}(K,\eta)$.
	
	\begin{theorem}\label{thm-4-1}
		Suppose $m\in\{1,\ldots,n-1\}$, $0<q\le n$. Let $L$ be a fixed convex body containing the origin in its interior,  and let $\mu$ be a finite even Borel measure with $|\mu|>0$.
		If $K_0\in\mathcal{K}_e^n$ is a maximizer for the problem:
		\[
		J^{(L)}_{m,q}(K_0)=\sup\{J^{(L)}_{m,q}(K): K\in\mathcal{K}_e^n\} ,
		\]
		then there exists a convex body $\widetilde K\in\mathcal{K}_e^n$ such that
		\[
		\mu(\cdot)=\widetilde{A}^{(L)}_{m,q}(\widetilde K,\cdot).
		\]
	\end{theorem}
	
	\begin{proof}
		Suppose $f\in C(S^{n-1})$ is even. Let $h_0=h_{K_0}$ and $\Omega=S^{n-1}$, and define $h_t$ by
		\[
		\log h_t=\log h_0+t f,\qquad t\in\mathbb{R}.
		\]
		Define the logarithmic family
		\[
		K_t=\{x\in\mathbb{R}^n : x\cdot v\le h_t(v)\ \text{for all } v\in S^{n-1}\},\qquad t\in\mathbb{R},
		\]
		and observe that $h_{K_t}\le h_t$.
		Consider the function
		\[
		\phi(t)=\frac{1}{m q}\log \widetilde{\Psi}_{m,q}(K_t,L)-\frac1{|\mu|}\int_{S^{n-1}}\log h_t(v)d\mu(v),
		\qquad t\in\mathbb{R}.
		\]
		From the fact that $h_{K_t}\le h_t$ and the assumption that
		$J^{(L)}_{m,q}(K_0)=\sup\{J^{(L)}_{m,q}(K,): K\in\mathcal{K}_e^n\}$, we have
		\[
		\phi(0)=J^{(L)}_{m,q}(K_0)\ge J^{(L)}_{m,q}(K_t)\ge \phi(t),
		\]
		for all $t$. Theorem~\ref{thm-3-4} tells us that $\phi(t)$ is differentiable at $t=0$.
		From this fact, together with the definition of $J_{m,q}^{(L)}$, the variational formula for $K$, it follows that
		\[
		\begin{aligned}
			0 &= \left. \frac{d}{dt} \phi(t) \right|_{t=0} \\
			&= \left. \frac{d}{dt} \left( \frac{1}{mq}\log \widetilde{\Psi}_{m,q}(K_t,L)-\frac1{|\mu|}\int_{S^{n-1}}\log h_t(v)d\mu(v)\right) \right|_{t=0} \\
			&= -\frac{1}{|\mu|} \int_{S^{n-1}} f(v)d\mu(v) + \frac{1}{m \widetilde{\Psi}_{m,q}(K_0, L)} \int_{S^{n-1}} f(v) d\widetilde{A}_{m,q}^{(L)}(K_0, v) \\
			&= \frac{1}{|\mu|} \left( -\int_{S^{n-1}} f(v)  d\mu(v) + \frac{|\mu|}{m \widetilde{\Psi}_{m,q}(K_0, L)} \int_{S^{n-1}} f(v) d\widetilde{A}_{m,q}^{(L)}(K_0, v) \right).
		\end{aligned}
		\]
		Let $\widetilde K=cK_0$, with
		\[
		c=\left(\frac{|\mu|}{m\,\widetilde{\Psi}_{m,q}(K_0,L)}\right)^{\frac{1}{q m}}.
		\]
		Then, by the $(mq)$-homogeneity of $\widetilde{\Psi}_{m,q}(\cdot,L)$ and $\widetilde{A}_{m,q}^{(L)}$, we obtain
		\[
		\int_{S^{n-1}} f(u)d\widetilde{A}_{m,q}^{(L)}(\widetilde K ,u)
		= \int_{S^{n-1}}f(u)d\mu(u).
		\]
		Since $f\in C(S^{n-1})$ is even but arbitrary, and the measure $\mu$ is even, the desired result follows.
	\end{proof}

	\section{Solution to the Minkowski problem}
	
	\quad Our aim in this section is to prove Theorem~\ref{thm-1-1} and Theorem~\ref{thm-1-2} concerning the existence of solutions to the Minkowski problems for mixed centro-section measures
	$\widetilde{A}^{(L)}_{m,q}(K,\eta)$, where $m\in\{1,\ldots,n-1\}$, with the ranges of $q$ specified in the corresponding statements below. We will deal with this problem under the assumption that the given measure satisfies the subspace concentration inequality.
	\[
	\frac{\mu(\xi \cap S^{n-1})}{|\mu|} < \frac{\dim \xi}{n}
	,\]	each nonzero proper subspace $\xi$ of $\mathbb{R}^n$.
	
	By Theorem~\ref{thm-4-1} , the solution of the Minkowski problem for the mixed centro-section measures $\widetilde{A}^{(L)}_{m,q}(K,\cdot)$ can be reduced to finding the maximizers of \eqref{4-1}. To achieve this, we require an estimate of the entropy $E_\mu(\cdot)$ in \cite{2019-Boroczky-AdvMath} (see \cite{2013-Boroczky-JAMS} for an earlier version), and thus we omit the proof.
	
	\begin{lemma}[\textnormal{\cite[Lemma 4.2]{2019-Boroczky-AdvMath}}]\label{lem-5-1}
		Let $\varepsilon_0>0$.
		Suppose that $(e_{1,l},\ldots,e_{n,l})$, where $l=1,2,\ldots$, is a sequence of ordered orthonormal bases of $\mathbb{R}^n$
		converging to the ordered orthonormal basis $(e_1,\ldots,e_n)$, while $(a_{1,l},\ldots,a_{n,l})$ is a sequence of $n$-tuples satisfying, for all $l$,
		\[
		0<a_{1,l}\le a_{2,l}\le \cdots \le a_{n,l}\quad\text{and}\quad a_{n,l}>\varepsilon_0.
		\]
		For each $l=1,2,\ldots$, let
		\[
		Q_l=\left\{x\in\mathbb{R}^n:\frac{|x\cdot e_{1,l}|^2}{a_{1,l}^2}+\cdots+\frac{|x\cdot e_{n,l}|^2}{a_{n,l}^2}\le 1\right\}
		\]
		denote the ellipsoid generated by the $(e_{1,l},\ldots,e_{n,l})$ and $(a_{1,l},\ldots,a_{n,l})$.
		Let $\mu$ be a Borel measure on $S^{n-1}$ that satisfies the strict subspace concentration inequality
		\[
		\frac{\mu(\xi\cap S^{n-1})}{|\mu|}<\frac{\dim\xi}{n},
		\]
		for each nonzero proper subspace $\xi$.
		Then there exists $t_0,\delta_0,l_0>0$ and $c_{n,\varepsilon_0,t_0,\delta_0}$ independent of $l$, such that, for each $l>l_0$,
		\[
		E_\mu(Q_l)\le -\frac{\log(a_{1,l}\cdots a_{n,l})}{n}+t_0\log a_{1,l}+c_{n,\varepsilon_0,t_0,\delta_0}.
		\]
	\end{lemma}
	\begin{theorem}\label{thm-5-2}
		Let $L$ be a fixed convex body containing the origin in its interior. Let $q\in(0,n]$ and $\mu$ be an even finite Borel measure on $S^{n-1}$ with $|\mu|>0$.
		If $\mu$ satisfies the subspace mass inequality
		\begin{equation*}
			\frac{\mu(\xi\cap S^{n-1})}{|\mu|}
			<\frac{\dim\xi}{n},
			\qquad \text{for each nonzero proper subspace }\xi\text{ of }\mathbb{R}^n,
		\end{equation*}
		then there exists a $K_q\in\mathcal{K}_e^n$ such that
		\[
		J^{(L)}_{m,q}(K_q)=\sup\{J^{(L)}_{m,q}(K): K\in\mathcal{K}_e^n\}.
		\]
	\end{theorem}
	
	\begin{proof}
		Fix $q\in(0,n]$, let $(K_j)$ be a sequence in $\mathcal{K}_e^n$ such that
		\[
		\lim_{j\to\infty}J^{(L)}_{m,q}(K_j)=\sup\{J^{(L)}_{m,q}(K): K\in\mathcal{K}_e^n\}.
		\]
		Since $J^{(L)}_{m,q}$ is $0$-homogeneous, we may assume the diameter of $K_j$ is $1$, for all $j$.
		By the Blaschke selection theorem, there is a subsequence of $(K_j)$, still denoted by $(K_j)$ for convenience,
		such that $K_j$ converges to some origin-symmetric compact convex set $K_q$.
		
		It remains to prove  that $K_q$ contains the origin in its interior, i.e., $K_q\in\mathcal{K}_e^n$.
		We argue by contradiction and assume  $K_q$ is contained in a proper subspace of $\mathbb{R}^n$.
		Since $K_j\in\mathcal{K}_e^n$, from John's theorem, there exists an origin-centered ellipsoid $Q_j$ such that
		\begin{equation}\label{5-1}
			Q_j\subset K_j\subset \sqrt{n}\,Q_j,\qquad j=1,2,\ldots.
		\end{equation}
		By the principal axis theorem, each ellipsoid $Q_j$ can be written as
		\[
		Q_j=\left\{x\in\mathbb{R}^n:\frac{|x\cdot e_{1,j}|^2}{a_{1,j}^2}+\frac{|x\cdot e_{2,j}|^2}{a_{2,j}^2}+\cdots+\frac{|x\cdot e_{n,j}|^2}{a_{n,j}^2}\le 1\right\},
		\]
		where $(e_{1,j},\ldots,e_{n,j})$ is an ordered orthonormal basis of $\mathbb{R}^n$, and
		\[
		0<a_{1,j}\le a_{2,j}\le \cdots \le a_{n,j}\le 1.
		\]
		The fact that $a_{n,j}\le 1$ follows from the fact that $Q_j\subset K_j$ together with $\operatorname{diam}(K_j)=1$. Since $(e_{1,j},\ldots,e_{n,j})$ is an orthonormal basis of $\mathbb{R}^n$, there exists a subsequence of
		$(e_{1,j},\ldots,e_{n,j})$ converging to an orthonormal basis of $\mathbb{R}^n$ that we denote by $(e_1,\ldots,e_n)$.
		For each $i=1,\ldots,n$, each sequence $(a_{i,j})$ is bounded and thus has a convergent subsequence.
		By successively passing to subsequences, we may assume that the corresponding ordered orthonormal bases $(e_{1,j},\ldots,e_{n,j})$ converge to an ordered orthonormal basis $({e_1,\ldots,e_n})$ of $\mathbb{R}^n$, with $a_{i,j}$ converging to say $a_i$,
		for every $i=1,\ldots,n$.
		
		Since we have assumed that $K_q$ is contained in some lower dimensional subspace and $a_{i,j}$ converges to $a_i$,
		there is an integer $k\in\{1,\ldots,n-1\}$ satisfying
		\[
		a_1=a_2=\cdots=a_k=0,\qquad \text{while}\qquad a_{k+1},\ldots,a_n>0.
		\]
		
		Since $K_j\subset \sqrt{n}Q_j$ and $\operatorname{diam}K_j=1$, we have $a_{n,j}\ge \frac{1}{2\sqrt{n}}$.
		According to Lemma~\ref{lem-5-1}, there exist positive real numbers $t_0,\delta_0,j_0$ and $c_{n,t_0,\delta_0}$ such that, for $j>j_0$,
		\begin{equation}\label{5-2}
			E_\mu(Q_j)\le -\frac{1}{n}\log(a_{1,j}\cdots a_{n,j})+t_0\log a_{1,j}+c_{n,t_0,\delta_0},
		\end{equation}
		where $c_{n,t_0,\delta_0}$ is independent of $j$.
		Since \(L\in\mathcal{K}_o^n\) is fixed, the continuity of the section-volume function
		$
		\xi\mapsto \operatorname{vol}_m(L\cap\xi)
		$
		on the compact space \(G(n,m)\) implies that
		$ C_L:=\max_{\xi\in G(n,m)}\operatorname{vol}_m(L\cap\xi)<\infty.$
		Consequently,
		$
		\operatorname{vol}_m(L\cap\xi)\le C_L
		$
		for every \(\xi\in G(n,m)\).
		
		Recall that $\tilde{\Psi}_{m,n}(K, L)$ is $\mathrm{SL}(n)$ invariant.
		This and the $mn$-homogeneity of $\widetilde{\Psi}_{m,n}(Q_j)$ yield
		$$
		\widetilde{\Psi}_{m,n}(Q_j,L)=\widetilde{\Psi}_{m,n}(Q_j),$$
		and
		\begin{equation}\label{5-3}
			\widetilde{\Psi}_{m,n}(Q_j,L)=(a_{1,j}\cdots a_{n,j})^m\,\widetilde{\Psi}_{m,n}(B^n),
		\end{equation}
		where $B^n$ denotes the unit ball. We obtain
		$$ \frac{1}{n m}\log\widetilde{\Psi}_{m,n}(Q_j,L)
		=\frac{\log(a_{1,j}\cdots a_{n,j})}{n}+\frac{1}{n m}\log\widetilde{\Psi}_{m,n}(B^n).$$
		If $q\in(0,n)$, we have $$
		\widetilde{\Psi}_{m,q}(Q_j,L)\leq C_L^{n-q}\widetilde{\Psi}_{m,q}(Q_j).$$ By \eqref{5-3} and the monotonicity of $L^p$ norm, we have
		\begin{equation}
			\begin{aligned}
				\frac{1}{q m}\log\widetilde{\Psi}_{m,q}(Q_j,L)
				&\le \frac{1}{q m}\log\bigl(C_L^{n-q}\widetilde{\Psi}_{m,q}(Q_j)\bigr) \\
				&\leq \frac{n-q}{mq}\log C_L + \frac{1}{n m}\log\widetilde{\Psi}_{m,n}(Q_j) \\
				&= \frac{n-q}{mq}\log C_L + \frac{\log(a_{1,j}\cdots a_{n,j})}{n} + \frac{1}{n m}\log\widetilde{\Psi}_{m,n}(B_n).
			\end{aligned}
			\label{5-4}
		\end{equation}
		Then, by \eqref{5-1}, \eqref{5-2} and \eqref{5-4}, for $j>j_0$, we have
		\[
		\begin{aligned}
			J^{(L)}_{m,q}(K_j)
			&\le \frac{n-q}{mq}\log{C_L}+\log\sqrt{n}+\frac{\log(a_{1,j}\cdots a_{n,j})}{n}
			+\frac{1}{n m}\log\widetilde{\Psi}_{m,n}(B_n)+E_\mu(Q_j)\\
			&\le t_0\log a_{1,j}+c_{n,t_0,\delta_0}+\log\sqrt{n}+\frac{n-q}{mq}\log{C_L}+\frac{1}{n m}\log\widetilde{\Psi}_{m,n}(B_n)\to -\infty,
		\end{aligned}\] \text{as} $j\to\infty$.
		
		However, since $J^{(L)}_{m,q}(K_j)$ converges to the supremum, we have
		\[
		-\infty=\lim_{j\to\infty}J^{(L)}_{m,q}(K_j)\ge J^{(L)}_{m,q}(B_n)>-\infty,
		\]
		which yields the contradiction. Thus $K_q$ has nonempty interior, that is $K_q\in\mathcal{K}_e^n.$
		For the endpoint $q = n$, the argument is even simpler because
		\[
		\widetilde{\Psi}_{m,n}(Q_j, L) = \widetilde{\Psi}_{m,n}(Q_j) = (a_{1,j} \cdots a_{n,j})^m \widetilde{\Psi}_{m,n}(B_n).
		\]
		Thus the same contradiction follows without introducing $C_L.$
		
		Since $K_j\to K_q$ in the Hausdorff metric and $K_q\in\mathcal K_e^n$,
		the continuity of $\widetilde{\Psi}_{m,q}(\,\cdot\,,L)$ and
		$E_\mu(\,\cdot\,)$ implies that $J_{m,q}^{(L)}$ attains its supremum at
		$K_q$. Hence, $K_q$ is a maximizer.
	\end{proof}
	
	\textbf{Proof of Theorem~\ref{thm-1-1}.} It follows immediately from Theorem~\ref{thm-5-2} and Theorem~\ref{thm-4-1}.\quad $\square$
	
	\vspace{0.2cm}
	
	\textbf{A refined subspace mass condition for $m=n-1$}
	
	In this subsection, for the case where $m = n-1$, we will solve the Minkowski problem for $\widetilde{A}_{n-1,q}^{(L)}(K,\cdot)$, under a refined and strictly weaker subspace mass condition. We say an even finite Borel measure $\mu$ satisfies the \textit{refined subspace mass inequality}, if for all $\xi_k \in G(n,k)$, where $k = 1,\dots,n-1$,
	\begin{equation}\label{5-5}
		\frac{\mu(\xi_k \cap S^{n-1})}{|\mu|} < \frac{k}{n-1} - \frac{(k + q - n)_+}{q(n-1)}, 
	\end{equation}
	where $x_+ = \max\{x,0\}$ denotes the positive part of $x\in\mathbb R$. Observing that when $0 < q < n$,
	\[
	\frac{k}{n} < \frac{k}{n-1} - \frac{(k + q - n)_+}{q(n-1)} \leq \frac{k}{n-1},
	\]
	\eqref{5-5} is strictly weaker than (1.6). In particular, when $q \in (0,1]$, \eqref{5-5} becomes
	\begin{equation}\label{5-6}
		\frac{\mu(\xi_k \cap S^{n-1})}{|\mu|} < \frac{k}{n-1},\quad \text{for all } \xi_k \in G(n,k). \end{equation}
	
	Recall that
	\begin{equation}\label{5-7}
		\widetilde{\Psi}_{n-1,q}(K,L) = \frac{1}{n\omega_n}\int_{S^{n-1}} \rho_{\mathrm{I}K}^q(u)\rho_{\mathrm{I}L}^{n-q}(u)\,du = \frac{1}{\omega_n}\widetilde{V}_{n-q}(\mathrm{I}K,\mathrm{I}L). \end{equation}
	This allows us to optimally estimate $\widetilde{\Psi}_{n-1,q}(Q,L)$ of an ellipsoid $Q$.

	\begin{lemma}\label{lem-5-3}{\cite[Lemma 6.3]{2026-Cai-AdvMath}}
		Let $q \in (l - 1, l)$, where $l \in [1, n]$ is an integer. Suppose $Q$ is an ellipsoid given by
		\[
		Q = \left\{ x \in \mathbb{R}^n : \frac{|x \cdot e_1|^2}{a_1^2} + \frac{|x \cdot e_2|^2}{a_2^2} + \cdots + \frac{|x \cdot e_n|^2}{a_n^2} \leq 1 \right\},
		\]
		where $0 < a_1 \leq a_2 \leq \cdots \leq a_n \leq 1$ and $e_1, \dots, e_n$ is an orthonormal basis. Then
		\[
		\widetilde{\Psi}_{n-1,q}(Q) \leq c_{q,l,n} \left(a_1 \cdots a_{n-l}\right)^q a_{n-l+1}^{l-1} \left(a_{n-l+2} \cdots a_n\right)^{q-1}, 
		\]
		where
		\[
		c_{q,l,n}=
		\begin{cases}
			\displaystyle \frac{2q\omega_{n-1}^{q+1}}{n(q - n + 1)\omega_n} & \text{when } l = n,\\[4pt]
			\displaystyle \frac{2^{\frac{n-q}{2}+1} q(n - l)\omega_{n-1}^q\omega_{n-l}\omega_{l-1}}{n(l - q)(q - l + 1)\omega_n} & \text{when } l \leq n-1.
		\end{cases}\]
	\end{lemma}
	Set
	\[
	r_L := \min_{u\in S^{n-1}} \rho_{\mathrm{I}L}(u) > 0, \qquad M_L := \max_{u\in S^{n-1}} \rho_{\mathrm{I}L}(u) < \infty.
	\]
	Since \( L \in \mathcal{K}_o^n \) is fixed, both \( r_L \) and \( M_L \) are positive constants depending only on \( L \).
	\begin{lemma}\label{lem-5-4}
		Let \( L \in \mathcal{K}_o^n \) be fixed and let \( 0 < q < n \). Then, for every \( K \in \mathcal{K}_o^n \),
		\[
		r_L^{n-q} \widetilde{\Psi}_{n-1,q}(K) \leq \widetilde{\Psi}_{n-1,q}(K,L) \leq M_L^{n-q} \widetilde{\Psi}_{n-1,q}(K).
		\]
	\end{lemma}
	
	\begin{proof}
		
		By (5.7) and  the definitions of \( r_L \) and \( M_L \), we have
		
		\[
		r_L^{n-q} \frac{1}{n\omega_n} \int_{S^{n-1}} \rho_{\mathrm{I}K}(u)^q du \leq \widetilde{\Psi}_{n-1,q}(K,L) \leq M_L^{n-q} \frac{1}{n\omega_n} \int_{S^{n-1}} \rho_{\mathrm{I}K}(u)^q du.
		\]
		Since
		\[
		\widetilde{\Psi}_{n-1,q}(K) = \frac{1}{n\omega_n} \int_{S^{n-1}} \rho_{\mathrm{I}K}(u)^q du,
		\]
		the desired inequalities follow.
	\end{proof}
	
	Combining Lemma~\ref{lem-5-3} with the upper estimate in Lemma~\ref{lem-5-4}, we obtain the following estimate for the mixed centro-section integral of an ellipsoid.
	\begin{lemma}\label{lem-5-5}
		Let $q \in (l - 1, l)$, where $l \in \{1,...,n\}.$  Let $L \in \mathcal{K}_o^n$ be fixed. Suppose $Q$ is an ellipsoid given by
		\[
		Q = \left\{ x \in \mathbb{R}^n : \frac{|x \cdot e_1|^2}{a_1^2} + \frac{|x \cdot e_2|^2}{a_2^2} + \cdots + \frac{|x \cdot e_n|^2}{a_n^2} \leq 1 \right\},
		\]
		where 
		$$0 < a_1 \leq a_2 \leq \cdots \leq a_n \leq 1$$ and $e_1, \dots, e_n$ is an orthonormal basis. Then
		\[
		\widetilde{\Psi}_{n-1,q}(Q, L) \leq c_{q,l,n,L} \left(a_1 \cdots a_{n-l}\right)^q a_{n-l+1}^{l-1} \left(a_{n-l+2} \cdots a_n\right)^{q-1},
		\]
		where
		\begin{equation}\label{5-8}
			c_{q,l,n,L}=
			\begin{cases}
				\displaystyle \frac{2q\omega_{n-1}^{q+1}}{n(q - n + 1)\omega_n} M_L^{n-q}, & l = n,\\[4pt]
				\displaystyle \frac{2^{\frac{n-q}{2}+1} q(n - l)\omega_{n-1}^q\omega_{n-l}\omega_{l-1}}{n(l - q)(q - l + 1)\omega_n} M_L^{n-q}, & l \leq n - 1.
			\end{cases}
		\end{equation}
		Clearly, the term $\left(a_1 \cdots a_{n-l}\right)^q$ in the above inequality will disappear if $l = n$, and the term $\left(a_{n-l+2} \cdots a_n\right)^{q-1}$ will disappear if $l = 1$.
	\end{lemma}
	
	The proof also requires the following entropy estimate.
	\begin{lemma}\label{lem-5-6}{\cite[Lemma 6.6]{2026-Cai-AdvMath}}
		Suppose $0<q<n$ and $\mu$ is a non-zero finite Borel measure on $S^{n-1}$ satisfying the refined subspace mass inequality
		$$
		\frac{\mu(\xi_k\cap S^{n-1})}{|\mu|} < \frac{k}{n-1} - \frac{(k+q-n)_+}{q(n-1)},\quad \text{for all }\xi_k\in G(n,k),
		$$
		where $k\in\{1,\dots,n-1\}$. Moreover, suppose $(e_{1,j},\dots,e_{n,j})$, for $j=1,2,\dots$, is a sequence of ordered orthonormal bases in $\mathbb R^n$ converging to the orthonormal basis $(e_1,\dots,e_n)$. Let $\varepsilon_0>0$ and a sequence of ellipsoids
		\[
		Q_j = \left\{x\in\mathbb R^n: \frac{|x\cdot e_{1,j}|^2}{a_{1,j}^2}+\cdots+\frac{|x\cdot e_{n,j}|^2}{a_{n,j}^2}\leq 1\right\},\quad j=1,2,\dots,
		\]
		with
		\[
		0 < a_{1,j}\leq a_{2,j}\leq \cdots\leq a_{n,j},\quad \text{and}\quad a_{n,j}\geq \varepsilon_0.
		\]
		Then, there are $j_0,\delta_0,t_0>0$ and a constant $c_{\varepsilon_0,\delta_0,t_0}$ independent of $j$, such that for all $j>j_0$
		\begin{equation}\label{5-9}
			\begin{aligned}
				E_\mu(Q_j) \leq& -\frac{1}{n-1}\log(a_{1,j}\cdots a_{n-\lfloor q\rfloor-1,j}) \\
				&- \frac{\lfloor q\rfloor}{q(n-1)}\log a_{n-\lfloor q\rfloor,j} - \frac{q-1}{q(n-1)}\log(a_{n-\lfloor q\rfloor+1,j}\cdots a_{n,j}) \\
				&+ t_0\log a_{1,j} + c_{\varepsilon_0,\delta_0,t_0}.
			\end{aligned}
		\end{equation}
	\end{lemma}
	
	The first term of \eqref{5-9} should be understood as disappearing if $\lfloor q\rfloor=n-1$, and the third term of \eqref{5-9} should also be understood as disappearing if $\lfloor q\rfloor=0$.
	\begin{theorem}\label{thm-5-7}
		Let \( L \in \mathcal{K}_o^n \) be fixed. Suppose \( q \in (0,n) \), and let \( \mu \) be a non-zero finite even Borel measure on \( S^{n-1} \). Assume that for every \( \xi_k \in G(n,k) \), \( k = 1,\dots,n-1 \),
		\[
		\frac{\mu(\xi_k \cap S^{n-1})}{|\mu|} < \frac{k}{n-1} - \frac{(k+q-n)_+}{q(n-1)}.
		\]
		Then, there exists \( K_q \in \mathcal{K}_e^n \) such that
		\[
		J_{n-1,q}^{(L)}(K_q) = \sup\bigl\{ J_{n-1,q}^{(L)}(K) : K \in \mathcal{K}_e^n \bigr\}.
		\]
	\end{theorem}
	\begin{proof} Let $(K_j)$ be a sequence in $\mathcal K_e^n$ such that
		\[
		\lim_{j\to\infty}J_{n-1,q}^{(L)}(K_j)
		=
		\sup\big\{J_{n-1,q}^{(L)}(K):K\in\mathcal K_e^n\big\}.
		\]
		Since $J_{n-1,q}^{(L)}$ is $0$-homogeneous, we may assume that
		$\operatorname{diam}K_j=1$ for all $j$. By the Blaschke selection theorem, after passing
		to a subsequence, $(K_j)$ converges to an origin-symmetric compact convex set $K_q$.
		It remains to prove that $K_q\in\mathcal K_e^n$. Suppose to the contrary that $K_q$
		has empty interior and is contained in a proper subspace of $\mathbb R^n$. By John's
		theorem, for each $j$ there is an origin-centered ellipsoid $Q_j$ such that
		\begin{equation}\label{5-10}
			Q_j\subset K_j\subset \sqrt n\,Q_j.
		\end{equation}
		For each $j$, there is an orthonormal basis $e_{1,j},\ldots,e_{n,j}$ such that
		\[
		Q_j=
		\left\{
		x\in\mathbb R^n:
		\frac{|x\cdot e_{1,j}|^2}{a_{1,j}^2}
		+\cdots+
		\frac{|x\cdot e_{n,j}|^2}{a_{n,j}^2}
		\le 1
		\right\},
		\]
		where
		\[
		0<a_{1,j}\le a_{2,j}\le\cdots\le a_{n,j}\le1,
		\qquad
		a_{n,j}\ge \frac{1}{2\sqrt n}.
		\]
		Passing to further subsequences, we may assume that the ordered orthonormal bases and
		all the semiaxis sequences converge. Since $K_q$ is contained in a proper subspace,
		necessarily $a_{1,j}\to0$ as $j\to\infty$.
		
		Set $l=\lfloor q\rfloor+1$. Then $l\in\{1,\ldots,n\}$ and $q\in[l-1,l)$. Since
		$Q_j\subset K_j$, we have $E_\mu(K_j)\le E_\mu(Q_j)$. Lemma~\ref{lem-5-6}, with
		$\varepsilon_0=1/(2\sqrt n)$, therefore yields constants $t_0>0$ and $C_0$,
		independent of $j$, such that
		\begin{equation}\label{5-11}
			\begin{aligned}
				E_\mu(K_j)
				&\le E_\mu(Q_j) \\
				&\le
				-\frac{1}{n-1}
				\log\!\bigl(a_{1,j}\cdots a_{n-l,j}\bigr)
				-\frac{l-1}{q(n-1)}
				\log a_{n-l+1,j} \\
				&\quad
				-\frac{q-1}{q(n-1)}
				\log\!\bigl(a_{n-l+2,j}\cdots a_{n,j}\bigr)
				+t_0\log a_{1,j}+c_{\varepsilon_0,\delta_0,t_0} .
			\end{aligned}
		\end{equation}
		Choose
		$q_0\in(q,l)\subset(l-1,l)$ sufficiently close to $q$ so that
		\begin{equation}\label{5-12}
			\frac{l-1}{n-1}\left(\frac1{q_0}-\frac1q\right)+t_0>0.
		\end{equation}
		In particular, $q_0$ is not an integer.
		Applying H\"older's inequality directly to \eqref{5-7}, we obtain
		\begin{align*}
			\widetilde\Psi_{n-1,q}(Q_j,L)
			&=
			\frac1{n\omega_n}
			\int_{S^{n-1}}
			\big(\rho_{\mathrm{I}Q_j}^{q_0}\rho_{\mathrm{I}L}^{n-q_0}\big)^{q/q_0}
			\big(\rho_{\mathrm{I}L}^{n}\big)^{1-q/q_0}\,du\\
			&\le
			\widetilde\Psi_{n-1,q_0}(Q_j,L)^{q/q_0}
			\left(
			\frac1{n\omega_n}
			\int_{S^{n-1}}\rho_{\mathrm{I}L}(u)^n\,du
			\right)^{1-q/q_0}\\
			&=
			\widetilde\Psi_{n-1,q_0}(Q_j,L)^{q/q_0}
			\left(\frac{V(\mathrm{I}L)}{\omega_n}\right)^{1-q/q_0}.
		\end{align*}
		Consequently,
		\begin{equation}\label{5-13}
			\frac1q\log\widetilde\Psi_{n-1,q}(Q_j,L)
			\le
			\frac1{q_0}\log\widetilde\Psi_{n-1,q_0}(Q_j,L)
			+
			\left(\frac1q-\frac1{q_0}\right)
			\log\frac{V(\mathrm{I}L)}{\omega_n}.
		\end{equation}
		Since $q_0\in(l-1,l)$,  Lemma~\ref{lem-5-5} applies to
		$\widetilde\Psi_{n-1,q_0}(Q_j,L)$. Combining Lemma~\ref{lem-5-5} with \eqref{5-13}, we obtain
		\begin{equation}\label{5-14}
			\begin{aligned}
				\frac{1}{q}\log \widetilde\Psi_{n-1,q}(Q_j,L)
				&\le
				\log\!\bigl(a_{1,j}\cdots a_{n-l,j}\bigr)
				+\frac{l-1}{q_0}\log a_{n-l+1,j} \\
				&\quad+
				\frac{q_0-1}{q_0}
				\log\!\bigl(a_{n-l+2,j}\cdots a_{n,j}\bigr)
				+\frac1{q_0}\log c_{q_0,l,n,L}
				+
				\left(\frac1q-\frac1{q_0}\right)
				\log\frac{V(\mathrm{I}L)}{\omega_n} .
			\end{aligned}
		\end{equation}
		By \eqref{5-10}, the monotonicity and homogeneity of
		$\widetilde\Psi_{n-1,q}(\cdot,L)$, combining \eqref{5-11} and \eqref{5-14}, we obtain
		\begin{equation}\label{5-15}
			\begin{aligned}
				J_{n-1,q}^{(L)}(K_j) &\leq \frac{1}{q(n-1)} \log \widetilde{\Psi}_{n-1,q}\big(\sqrt{n}\, Q_j, L\big) + E_\mu(Q_j) \\
				&= \frac{1}{q(n-1)} \log \widetilde{\Psi}_{n-1,q}(Q_j, L) + E_\mu(Q_j) + \frac{1}{2}\log n\\
				&\le
				\frac{l-1}{n-1}
				\left(\frac{1}{q_0}-\frac{1}{q}\right)
				\log a_{n-l+1,j} +
				\frac{1}{n-1}
				\left(\frac{1}{q}-\frac{1}{q_0}\right)
				\log\!\bigl(a_{n-l+2,j}\cdots a_{n,j}\bigr)	+t_0\log a_{1,j} \\
				&+\frac{1}{n-1}\left(\frac1{q_0}\log c_{q_0,l,n,L}
				+
				\left(\frac1q-\frac1{q_0}\right)
				\log\frac{V(\mathrm{I}L)}{\omega_n}\right)+c_{\varepsilon_0,\delta_0,t_0} +\frac{1}{2}\log n .
			\end{aligned}
		\end{equation}
		Since $q_0>q$ and
		$0<a_{i,j}\le1$, the second term on the right-hand side of \eqref{5-15} is non-positive.
		Moreover, $
		\frac{l-1}{n-1}
		\left(\frac1{q_0}-\frac1q\right)\le0,
		\;
		\log a_{1,j}\le\log a_{n-l+1,j}.
		$
		It follows that
		\[
		\begin{aligned}
			J_{n-1,q}^{(L)}(K_j)
			&\le{}
			\left[
			\frac{l-1}{n-1}
			\left(\frac{1}{q_0}-\frac{1}{q}\right)+t_0
			\right]\log a_{1,j} \\
			&+
			\frac{1}{n-1}
			\left[
			\frac{1}{q_0}\log c_{q_0,l,n,L}
			+
			\left(\frac{1}{q}-\frac{1}{q_0}\right)
			\log\frac{V(\mathrm{I}L)}{\omega_n}
			\right]
			+c_{\varepsilon_0,\delta_0,t_0}
			+\frac{1}{2}\log n .
		\end{aligned}
		\]
		By \eqref{5-12}, the coefficient in brackets is positive. Since $a_{1,j}\to0$, we conclude that
		\[
		\lim_{j\to\infty}J_{n-1,q}^{(L)}(K_j)=-\infty,
		\]
		which contradicts
		\[
		\lim_{j\to\infty}J_{n-1,q}^{(L)}(K_j)
		=
		\sup\big\{J_{n-1,q}^{(L)}(K):K\in\mathcal K_e^n\big\}
		\ge J_{n-1,q}^{(L)}(B^n)>-\infty.
		\]
		Thus $K_q\in\mathcal K_e^n$. Since $K_j\to K_q$ in the Hausdorff metric,
		the continuity of $\widetilde\Psi_{n-1,q}(\cdot,L)$ and $E_\mu(\cdot)$
		implies that $J_{n-1,q}^{(L)}$ attains its supremum at $K_q$.
		Hence $K_q$ is a maximizer.
	\end{proof}
	
	\noindent\textbf{Proof of Theorem 1.2.} It follows immediately from Theorem~\ref{thm-5-7} and Theorem~\ref{thm-4-1}. \quad $\square$
	
	\section{\texorpdfstring{SL($n$)-invariance and isoperimetric inequalities for the operator  $\mathbb{I}_{m,q}$}{I\_\{m,q\}}}
	\quad Throughout this section, invariance and covariance are considered under the
	action of the special linear group $\mathrm{SL}(n)$. We first prove that the
	mixed centro-section functional $\widetilde{\Psi}_{m,q}(K,L)$ is invariant
	under simultaneous $\mathrm{SL}(n)$-transformations of its two arguments,
	whereas the mixed bi-dual $m$-intersection bodies $\mathbb{I}_{m,q}(K,L)$ of index $q$ satisfy the relative covariance law
	$
	\mathbb{I}_{m,q}(\phi K,L)
	=
	\phi\,\mathbb{I}_{m,q}(K,\phi^{-1}L),
	$
	or, equivalently,
	$
	\mathbb{I}_{m,q}(\phi K,\phi L)
	=
	\phi\,\mathbb{I}_{m,q}(K,L).
	$
	We then establish sharp affine isoperimetric inequalities for the operator $\mathbb{I}_{m,q}$. These inequalities extend the affine
	isoperimetric inequality for the bi-dual $m$-intersection body
	operator $\mathbb{I}_m$ obtained in \cite{2025-Cai-AdvMath}. The latter is
	recovered when $q=1$, $q=n$, or $K=L$.
	
	\begin{theorem}\label{thm-6-1} Let $K$ and $L$ be origin--symmetric convex bodies in $\mathbb{R}^n,$ let $m\in\{1,\ldots,n-1\}$ and $q \in \mathbb{R}$.  Then $\widetilde{\Psi}_{m,q}(gK,gL)=\widetilde{\Psi}_{m,q}(K,L)$, provided $g \in\mathrm{SL}(n).$
	\end{theorem}
	\begin{proof}
		Let $\xi\in G(n,m)$ and $g\in \mathrm{SL}(n)$. Then the function
		\begin{equation}\label{6-1}
			\sigma_m(g,\xi)=\frac{\operatorname{vol}_m\bigl(g(M\cap \xi)\bigr)}
			{\operatorname{vol}_m(M\cap \xi)}
		\end{equation}
		is a multiplier function on the Grassmannian $G(n,m)$, where $M\in\mathcal{K}^n_o$.
		Let $g^{*}\nu_m$ denote the pull‑back of $\nu_m$ under the action
		$\xi\mapsto g\xi$ on $G(n,m)$; that is,
		\[
		g^{*}\nu_m(\omega)=\nu_m(g\omega)
		\]
		for every Borel set $\omega\subset G(n,m)$. The Radon–Nikodym derivative
		\begin{equation}\label{6-2}
			\sigma_{G(n,m)}(g,\xi)=
			\frac{d(g^{*}\nu_m)}{d\nu_m}(\xi)
		\end{equation}
		is also a multiplier on $G(n,m)$. A calculation presented in \cite{1971-Furstenberg-IsrJMath} yields
		\begin{equation}\label{6-3}
			\sigma_{G(n,m)}(g,\xi)=
			\sigma_m(g,\xi)^{-n}.
		\end{equation}
		Now, using \eqref{6-1}, \eqref{6-2}, and \eqref{6-3}, we obtain
		\[
		\begin{aligned}
			\widetilde{\Psi}_{m,q}(gK,gL)
			&=
			\int_{G(n,m)}
			\operatorname{vol}_m(gK\cap \xi)^q
			\operatorname{vol}_m(gL\cap \xi)^{n-q}
			\,d\nu_m(\xi) \\
			&=
			\int_{G(n,m)}
			\operatorname{vol}_m(gK\cap g\xi)^q
			\operatorname{vol}_m(gL\cap g\xi)^{n-q}
			\,d(g^{*}\nu_m)(\xi) \\
			&=
			\int_{G(n,m)}
			\bigl[
			\sigma_m(g,\xi)\operatorname{vol}_m(K\cap \xi)
			\bigr]^q
			\bigl[
			\sigma_m(g,\xi)\operatorname{vol}_m(L\cap \xi)
			\bigr]^{n-q} \\
			&\qquad\qquad\cdot
			\sigma_{G(n,m)}(g,\xi)
			\,d\nu_m(\xi) \\
			&=
			\int_{G(n,m)}
			\sigma_m(g,\xi)^n
			\operatorname{vol}_m(K\cap \xi)^q
			\operatorname{vol}_m(L\cap \xi)^{n-q} \\
			&\qquad\qquad\cdot
			\sigma_m(g,\xi)^{-n}
			\,d\nu_m(\xi) \\
			&=
			\int_{G(n,m)}
			\operatorname{vol}_m(K\cap \xi)^q
			\operatorname{vol}_m(L\cap \xi)^{n-q}
			\,d\nu_m(\xi) \\
			&=
			\widetilde{\Psi}_{m,q}(K,L).
		\end{aligned}
		\]
		This completes the proof.
	\end{proof}
The functional $\widetilde{\Psi}_{m,q}$ is invariant under simultaneous
$SL(n)$-transformations of its two arguments, but need not be invariant
when only one argument is transformed. The following example shows
that, even within the range $0<q<n$, there exist
$K_0,L\in\mathcal K_o^n$ and $g\in SL(n)$ such that
\begin{equation}\label{6-4}
	\widetilde{\Psi}_{m,q}(K_0,gL)
	\neq \widetilde{\Psi}_{m,q}(K_0,L).
\end{equation}
	
\begin{example}
Let
	$n=2$, $m=1$, $q=1$, and $K_0=B^2$. Take
	\[
	g=\operatorname{diag}(a,a^{-1})\in SL(2),
	\qquad a>1 .
	\]
	Then $gB^2$ is an origin-centered ellipse with
	$
	V(gB^2)=V(B^2).
	$
	For $\xi\in G(2,1)$, we have
	$
	\operatorname{vol}_1(gB^2\cap\xi)=2\rho_{gB^2}(\xi),
	$
	and therefore
	\[
	\widetilde{\Psi}_{1,1}(B^2,gB^2)
	=
	4\int_{G(2,1)}\rho_{gB^2}(\xi)d\nu_1(\xi).
	\]
	Since
	\[
	\int_{G(2,1)}\rho_{gB^2}(\xi)^2d\nu_1(\xi)=1
	\]
	and $\rho_{gB^2}$ is not constant, the strict Cauchy--Schwarz
	inequality gives
	\[
	\int_{G(2,1)}\rho_{gB^2}(\xi)d\nu_1(\xi)<1.
	\]
	Consequently,
	$	\widetilde{\Psi}_{1,1}(B^2,gB^2)
	<
	4
	=
	\widetilde{\Psi}_{1,1}(B^2,B^2).
	$
	Hence $\widetilde{\Psi}_{m,q}(K,L)$ is not affine invariant when
	only one argument is transformed.
\end{example}

	To study the properties of the operator $\mathbb{I}_{m,q}{(K,L)}$, we introduce the following definitions.
	
	For a fixed star body $K \in \mathcal{S}_o^n$, consider the mixed bi-dual $m$-intersection body $\mathbb{I}_{m,q}(K,L)$ with $ L \in \mathcal{S}_o^n.$
	We say that the family
	$
	\bigl\{\mathbb{I}_{m,q}(K,L) : L \in \mathcal{S}_o^n\bigr\}
	$
	is \textit{relatively affine covariant with respect to} $K$ if, for every $L \in \mathcal{S}_o^n$ and every $\varphi \in SL(n)$ such that
	\[
	\mathbb{I}_{m,q}(K,\varphi L) = \varphi \mathbb{I}_{m,q}(\varphi^{-1}K,L).
	\]
	Equivalently,
	\[
	\mathbb{I}_{m,q}(\varphi K,\varphi L) = \varphi \mathbb{I}_{m,q}(K,L).
	\]
	\begin{theorem}\label{thm-6-2}Let $K, L \in \mathcal{S}_o^n$, $\varphi \in SL(n)$, $m\in \{1,...,n-1\},$ and $q\in \mathbb{R}.$ Then the mixed bi-dual $m$-intersection bodies satisfy the following relative affine covariance
		\[
		\mathbb{I}_{m,q}(\varphi K,L) = \varphi \mathbb{I}_{m,q}{(K,\varphi^{-1}L)}.
		\]
		Equivalently,
		
		\[
		\mathbb{I}_{m,q}{(\varphi K,\varphi L)} = \varphi \mathbb{I}_{m,q}{(K,L)}.
		\]
	\end{theorem}
	
	\begin{proof}
		Let $ \varphi\in \mathrm{SL}(n),$  $K,L\in\mathcal{S}_o^n,$  and
		$M\in\mathcal{S}_o^n.$
		By the definition \eqref{1-4} of $\mathbb{I}_{m,q}(K,L)$ , using \eqref{2-5}, \eqref{6-1}, \eqref{6-2}, and \eqref{6-3},  we obtain
		\begin{align*}
			\widetilde{V}_m \big( \mathbb{I}_{m,q}(\varphi K,L), M \big)
			&= \frac{1}{n} \int_{S^{n-1}} \rho_M(u)^m \rho_{\mathbb{I}_{m,q}(\varphi K, L)}(u)^{n-m} du \\
			&= \frac{1}{n} \int_{S^{n-1}} \rho_M(u)^m \mathcal{R}_m^* \big( \operatorname{vol}_m(\varphi K \cap \cdot)^{q-1} \operatorname{vol}_m(L \cap \cdot)^{n-q} \big)(u) du \\
			&= \frac{m}{n} \int_{G(n,m)} \operatorname{vol}_m(M \cap \xi) \operatorname{vol}_m(\varphi K \cap \xi)^{q-1} \operatorname{vol}_m(L \cap \xi)^{n-q} d\nu_m(\xi) \\
			&= \frac{m}{n} \int_{G(n,m)} \operatorname{vol}_m(M \cap \varphi \xi) \operatorname{vol}_m(\varphi K \cap \varphi \xi)^{q-1} \operatorname{vol}_m(L \cap \varphi \xi)^{n-q} d(\varphi^{*}\nu_m)(\xi) \\
			&= \frac{m}{n} \int_{G(n,m)} \big[ \sigma_m(\varphi,\xi) \operatorname{vol}_m(\varphi^{-1}M \cap \xi) \big]
			\big[ \sigma_m(\varphi,\xi) \operatorname{vol}_m(K \cap \xi) \big]^{q-1} \\
			&\quad \cdot \big[ \sigma_m(\varphi,\xi) \operatorname{vol}_m(\varphi^{-1}L \cap \xi) \big]^{n-q} \sigma_m(\varphi,\xi)^{-n} d\nu_m(\xi) \\
			&= \frac{m}{n} \int_{G(n,m)} \operatorname{vol}_m(\varphi^{-1}M \cap \xi) \operatorname{vol}_m(K \cap \xi)^{q-1} \operatorname{vol}_m(\varphi^{-1}L \cap \xi)^{n-q} d\nu_m(\xi) \\
			&= \widetilde{V}_m \big( \mathbb{I}_{m,q}(K,\varphi^{-1}L), \varphi^{-1}M \big) \\
			&= \widetilde{V}_m \big( \varphi \mathbb{I}_{m,q}(K,\varphi^{-1}L), M \big).
		\end{align*}
		Now we have shown $\widetilde{V}_m\bigl(\mathbb{I}_{m,q}(\varphi K,L), M\bigr)=\widetilde{V}_m\bigl(\varphi \mathbb{I}_{m,q}(K,\varphi^{-1}L), M\bigr)$  for all $M\in\mathcal{S}_o^n$. Applying this and the dual mixed volume inequality \eqref{2-1} for $M=\mathbb{I}_{m,q}(\varphi K, L)$, we obtain
		\begin{align*}
			V\left(\mathbb{I}_{m,q}(\varphi K, L)\right)
			&= \widetilde{V}_m\left(\mathbb{I}_{m,q}(\varphi K,L), \mathbb{I}_{m,q}(\varphi K, L)\right) \\
			&= \widetilde{V}_m\left(\varphi \mathbb{I}_{m,q}(K,\varphi^{-1}L), \mathbb{I}_{m,q}(\varphi K,  L)\right)\\
			&\leq V\left(\varphi \mathbb{I}_{m,q}(K,\varphi^{-1}L)\right)^{(n-m)/n} V\left(\mathbb{I}_{m,q}(\varphi K,L)\right)^{m/n},
		\end{align*}
		which implies $	V\left(\mathbb{I}_{m,q}(\varphi K,L)\right) \leq V\left(\varphi \mathbb{I}_{m,q}(K,\varphi^{-1}L )\right)$.
		Similarly, taking $M =\varphi \mathbb{I}_{m,q}(K,\varphi^{-1}L)$, we have
		$V\left(\varphi \mathbb{I}_{m,q}(K,\varphi^{-1}L)\right) \leq V\left(\mathbb{I}_{m,q}(\varphi K,L)\right) $.
		Thus $V\left(\mathbb{I}_{m,q}(\varphi K,L)\right) = V\left(\varphi \mathbb{I}_{m,q}(K,\varphi^{-1}L)\right)$. Hence equality holds in \eqref{2-1}. We conclude that
		\[
		\mathbb{I}_{m,q}(\varphi K,L)
		= \varphi\,\mathbb{I}_{m,q}(K,\varphi^{-1}L).
		\]
		Equivalently,
		$$
		\mathbb{I}_{m,q}(\varphi K,\varphi L) = \varphi \mathbb{I}_{m,q}(K,L).$$
	\end{proof}
	In the  case $K = L$, the relative affine covariance in Theorem~\ref{thm-6-2} reduces to the ordinary affine covariance of the one-body bi-dual $m$-intersection body \cite{2025-Cai-AdvMath}. Thus, the following corollary is immediate.
	\begin{corollary}\label{cor-6-3}
		If $K = L \in \mathcal{S}_o^n$, then
		$
		\mathbb{I}_m(K) = \mathbb{I}_{m,q}(K,K)
		$
		and 
		$
		\mathbb{I}_m(\varphi K) = \varphi \mathbb{I}_m(K), \; \varphi \in SL(n),
		$
		where the bi-dual $m$-intersection operator $\mathbb{I}_m$ is defined by
		\[
		\rho_{\mathbb{I}_m K}(u)^{n-m} = \mathcal{R}_m^{*}\big(\operatorname{vol}_m(K \cap \cdot)^{n-1}\big)(u).
		\]
	\end{corollary}
	\begin{proof}
		By Theorem~\ref{thm-6-2}, in its equivalent form,
		\[
		\mathbb{I}_{m,q}(\varphi K,\varphi K) = \varphi \mathbb{I}_{m,q}(K,K).
		\]
		Therefore,
		\[
		\begin{aligned}
			\mathbb{I}_{m}(\varphi K)= \mathbb{I}_{m,q}(\varphi K,\varphi K) = \varphi \mathbb{I}_{m,q}(K,K)= \varphi \mathbb{I}_{m}(K).
		\end{aligned}
		\]
		This proves the affine covariance of the  operator $\mathbb{I}_{m}$. \end{proof}
	
	To establish the corresponding isoperimetric inequalities for the operator $\mathbb{I}_{m,q}(K,L)$, we need the following two lemmas.

	\begin{lemma}[\textnormal{\cite[Corollary 7.5]{2007-Gardner-AdvMath}}]\label{lem-6-4}
	Let $K\in\mathcal B^n$ satisfy $V(K)>0$, where $\mathcal B^n$
	denotes the class of bounded Borel sets in $\mathbb R^n$.
	For $1\le m\le n-1$, we have
		\[
		\int_{G(n,m)}\operatorname{vol}_m(K \cap \xi)^n d\nu_m(\xi) \leqslant \frac{\omega_m^n}{\omega_n^m} V(K)^m,
		\]
		with equality when $m > 1$ if and only if $K$ is an origin-symmetric $n$-dimensional ellipsoid, modulo a set of measure zero and when $m= 1$ if and only if $K$ is an origin-symmetric star-shaped set, modulo a set of measure zero.	
		
	\end{lemma}
	\begin{lemma}[\textnormal{\cite[Theorem 6.2]{2025-Cai-AdvMath}}]\label{lem-6-5}
		Let $K \in \mathcal{S}_o^n$,
		\(m\in\{2,\ldots,n-1\}\).  Then
		\[
		V \bigl(\mathbb{I}_m K\bigr) \leq \left( \frac{m \omega_m^n}{n \omega_n^m} \right)^{\frac{n}{n-m}} V(K)^{\frac{m(n-1)}{n-m}},
		\]
		with equality if and only if $K$ is an ellipsoid centered at the origin.
	\end{lemma}
	\noindent\textbf{Remark 6.6.}
		We emphasize that the result in Lemma~\ref{lem-6-5} is valid in the larger class of star bodies \(\mathcal{S}_o^n\), rather than only for convex bodies. This broader formulation will be used in the proof of Theorem~\ref{thm-6-8}. We also note that when $m=n$, equality holds in the inequality in Lemma~\ref{lem-6-4} for every bounded Borel set $K$.
	
	\begingroup
	\renewcommand{\thetheorem}{6.7}
	\begin{theorem}\label{thm-6-7}
		Let \(K,L\in\mathcal{S}_o^n\), let \(m=1\), and let
		\(1\leq q\leq n\). Then
		\begin{equation}\label{6-5}
			V\bigl(\mathbb{I}_{1,q}(K,L)\bigr)
			\leq
			\left(
			\frac{2^n}{n\omega_n}
			\right)^{\frac{n}{n-1}}
			V(K)^{\frac{q-1}{n-1}}
			V(L)^{\frac{n-q}{n-1}}.
		\end{equation}
		The equality cases are as follows:
		\begin{enumerate}
			\setlength{\topsep}{4pt}
			\setlength{\itemsep}{2pt}
			\setlength{\parsep}{0pt}
			\setlength{\parskip}{0pt}
			
			\item[(a)] If \(q=1\), equality holds if and only if \(L\) is an
			origin-symmetric star body, while \(K\) is arbitrary.
			
			\item[(b)] If \(q=n\), equality holds if and only if \(K\) is an
			origin-symmetric star body, while \(L\) is arbitrary.
			
			\item[(c)] If \(1<q<n\), equality holds if and only if \(K\) and \(L\)
			are origin-symmetric star bodies that are dilates of one another.
		\end{enumerate}
	\end{theorem}
	
	\begin{proof}
		For $A\in\mathcal{S}_o^n$ and $\xi\in G(n,1)$, set
		$ s_A(\xi)=\operatorname{vol}_1(A\cap\xi).$
		For $u\in S^{n-1}$, let $\xi_u=\operatorname{span}\{u\}$ and write $ s_A(u)=s_A(\xi_u).$
		Then
		$
		s_A(u)=\rho_A(u)+\rho_A(-u).$
		
		By the definition of the mixed bi-dual intersection body,
		\[
		\rho_{\mathbb{I}_{1,q}(K,L)}(u)^{n-1} = \frac{2}{n\omega_n} s_K(u)^{q-1} s_L(u)^{n-q}.
		\]
		Consequently,
		\[
		\begin{aligned}
			V\big(\mathbb{I}_{1,q}(K,L)\big) &= \frac{1}{n}\int_{S^{n-1}} \rho_{\mathbb{I}_{1,q}(K,L)}(u)^n du \\
			&= \frac{1}{n}\left(\frac{2}{n\omega_n}\right)^{\frac{n}{n-1}} \int_{S^{n-1}} s_K(u)^{\frac{n(q-1)}{n-1}} s_L(u)^{\frac{n(n-q)}{n-1}} du.
		\end{aligned}
		\]
		Suppose first that $1 < q < n$. Applying Hölder's inequality, we obtain
		\begin{equation}\label{6-6}
			V\big(\mathbb{I}_{1,q}(K,L)\big) \leq \frac{1}{n}\left(\frac{2}{n\omega_n}\right)^{\frac{n}{n-1}} \left(\int_{S^{n-1}} s_K(u)^n du\right)^{\frac{q-1}{n-1}} \left(\int_{S^{n-1}} s_L(u)^n du\right)^{\frac{n-q}{n-1}}.
		\end{equation}	
		By the convexity of the function $t\mapsto t^n$ on $[0,\infty)$, for every
		$u\in S^{n-1}$, we obtain
		\begin{equation}\label{6-7}
			\bigl(\rho_A(u)+\rho_A(-u)\bigr)^n
			\le
			2^{n-1}\bigl(\rho_A(u)^n+\rho_A(-u)^n\bigr),
		\end{equation}
		with equality if and only if $\rho_A(u)=\rho_A(-u)$.
		Integrating over $S^{n-1}$ and using $V(A) = \frac{1}{n}\int_{S^{n-1}} \rho_A(u)^n du$, we obtain
		\begin{equation}\label{6-8}
			\int_{S^{n-1}} s_A(u)^n du \leq 2^n n V(A).
		\end{equation}
		Substituting the above inequality into the previous estimate yields
		\[
		V\big(\mathbb{I}_{1,q}(K,L)\big) \leq \left(\frac{2^n}{n\omega_n}\right)^{\frac{n}{n-1}} V(K)^{\frac{q-1}{n-1}} V(L)^{\frac{n-q}{n-1}}.
		\]	
		We now determine the equality cases. Note that equality holds in \eqref{6-8} if and only if equality holds in \eqref{6-7} for almost every $u$, which, by the continuity of the radial function, is equivalent to $A$ being origin-symmetric.
		Let $1<q<n$. Equality holds if and only if equality occurs in Hölder’s inequality \eqref{6-6} and in \eqref{6-7}  for both $K$ and $L$. Hence  $K$ and $L$ are origin-symmetric and $s_K^n=cs_L^n$ for some $c>0.$ Since $s_K=2\rho_K$  and $s_L=2\rho_L$, the  bodies are dilates. The equality cases for $q=1$ and $q=n$ are immediate. The converse follows directly.
	\end{proof}

	\begingroup
\renewcommand{\thetheorem}{6.8}
\begin{theorem}\label{thm-6-8}
		Let \(K,L\in\mathcal{S}_o^n\), let
		\(m\in\{2,\ldots,n-1\}\), and let \(1\leq q\leq n\). Then
		\[
		V\bigl(\mathbb{I}_{m,q}(K,L)\bigr)
		\leq
		\left(
		\frac{m\omega_m^n}{n\omega_n^m}
		\right)^{\frac{n}{n-m}}
		V(K)^{\frac{m(q-1)}{n-m}}
		V(L)^{\frac{m(n-q)}{n-m}}.
		\]
		The equality cases are as follows:
		\begin{enumerate}
			\setlength{\topsep}{4pt}
			\setlength{\itemsep}{2pt}
			\setlength{\parsep}{0pt}
			\setlength{\parskip}{0pt}
			
			\item[(a)] If \(q=1\), equality holds if and only if \(L\) is an
			origin-centered ellipsoid, while \(K\) is arbitrary.
			
			\item[(b)] If \(q=n\), equality holds if and only if \(K\) is an
			origin-centered ellipsoid, while \(L\) is arbitrary.
			
			\item[(c)] If \(1<q<n\), equality holds if and only if \(K\) and \(L\)
			are origin-centered ellipsoids that are dilates of one another.
		\end{enumerate}
	\end{theorem}
	
	\begin{proof}
		For $A\in\mathcal{S}_o^n$ and $\xi\in G(n,m)$, set
		$$
		s_A(\xi) = \operatorname{vol}_m(A\cap\xi).
		$$
		Let $M\in\mathcal{S}_o^n$. By the duality of the Radon transform and the definition \eqref{1-4} of $\mathbb{I}_{m,q}(K,L)$, we have
		\[
		\widetilde V_m\bigl(\mathbb{I}_{m,q}(K,L),M\bigr)
		=
		\frac{m}{n}\int_{G(n,m)} s_M(\xi)s_K(\xi)^{q-1}s_L(\xi)^{n-q}\,d\nu_m(\xi).
		\]
		Applying the generalized Hölder inequality,  we obtain
		\[
		\begin{aligned}
			&\int_{G(n,m)} s_M s_K^{q-1}s_L^{n-q}\,d\nu_m\\
			&\leq \left(\int_{G(n,m)} s_M^n\,d\nu_m\right)^{1/n}
			\left(\int_{G(n,m)} s_K^n\,d\nu_m\right)^{(q-1)/n}
			\left(\int_{G(n,m)} s_L^n\,d\nu_m\right)^{(n-q)/n}.
		\end{aligned}
		\]
		By Lemma~\ref{lem-6-4}, for each $A\in\mathcal{S}_o^n$,
		\[
		\int_{G(n,m)} s_A(\xi)^n\,d\nu_m(\xi) \leq \frac{\omega_m^n}{\omega_n^m}V(A)^m.
		\]
		Consequently,
		\[
		\widetilde V_m\bigl(\mathbb{I}_{m,q}(K,L),M\bigr)
		\leq
		\frac{m\omega_m^n}{n\omega_n^m}V(M)^{m/n}V(K)^{m(q-1)/n}V(L)^{m(n-q)/n}.
		\]
		Taking $M = \mathbb{I}_{m,q}(K,L)$ and using
		$
		\widetilde V_m\bigl(\mathbb{I}_{m,q}(K,L),\mathbb{I}_{m,q}(K,L)\bigr) = V\bigl(\mathbb{I}_{m,q}(K,L)\bigr),
		$
		we obtain
		\[
		V\bigl(\mathbb{I}_{m,q}(K,L)\bigr)^{\frac{n-m}{n}}
		\leq
		\frac{m\omega_m^n}{n\omega_n^m}V(K)^{\frac{m(q-1)}{n}}V(L)^{\frac{m(n-q)}{n}}.
		\]
		This proves the asserted inequality. We next establish the equality conditions.
		
		In case (a), $q=1$ and $\mathbb{I}_{m,1}(K,L)=\mathbb{I}_{m}(L)$. Lemma~\ref{lem-6-5} gives the stated equality condition on $L$.
		
	In case (b), $q=n$ and $\mathbb{I}_{m,n}(K,L)=\mathbb{I}_{m}(K)$. Lemma~\ref{lem-6-5} gives the stated equality condition on $K$.
		
		For case (c), let $1<q<n$ and set $M=\mathbb{I}_{m,q}(K,L)$. The equality case forces equality in H\"older's inequality and in the inequalities involving $M$, $K$, and $L$ stated in Lemma~\ref{lem-6-4}. Hence all three bodies are origin-centered ellipsoids. Moreover,
		\[
		s_K(\xi)^n=c\,s_L(\xi)^n
		\]
		for some $c>0$ and $\nu_m$-almost every $\xi\in G(n,m)$. Since the functions $s_K$ and $s_L$ are continuous on
		$G(n,m)$,
		$
		s_A=\frac{1}{m}\mathcal{R}_m(\rho_A^m)
		$
		and $\mathcal{R}_m$ is injective on the space of even continuous functions \cite{2024-RB-JAM}, it follows that
		\[
		\rho_K=c^{1/(mn)}\rho_L.
		\]
		Thus $K$ and $L$ are dilates. Conversely, if $K$ and $L$ are origin-centered ellipsoids that are dilates of one another, then equality holds in all the preceding inequalities. This proves (c).	
	\end{proof}
\renewcommand{\thetheorem}{6.9}
\begin{corollary}
	Let $K\in\mathcal S_o^n$, $m\in\{1,\ldots,n-1\}$, and
	$1<q<n$. Define the operator $\mathbb I_{m,q}$ by
	\[
	\rho_{\mathbb I_{m,q}K}(u)^{n-m}
	=
	\mathcal R_m^*
	\bigl(\operatorname{vol}_m(K\cap\cdot)^{q-1}\bigr)(u),\quad \text{for}\; u\in S^{n-1}.
	\]
	Then
	\[
	\frac{V(\mathbb I_{m,q}K)}
	{V(\mathbb I_{m,q}B^n)}
	\le
	\left(\frac{V(K)}{\omega_n}\right)^{
		\frac{m(q-1)}{n-m}},
	\]
	with equality if and only if $K$ is a Euclidean ball
	centered at the origin.
	\end{corollary}

\begin{proof}
	Set $L_0=\omega_m^{-1/m}B^n$. Since
	$
	\operatorname{vol}_m(L_0\cap\xi)=1
	\;\text{for each}\; \xi\in G(n,m),
	$
	and
	$
	\mathbb I_{m,q}(K,L_0)=\mathbb I_{m,q}K.
	$
	Applying Theorem~\ref{thm-6-7}, Theorem~\ref{thm-6-8} and using
	$V(L_0)=\omega_n\omega_m^{-n/m}$, we obtain
	\[
	V(\mathbb I_{m,q}K)
	\le
	\omega_n
	\left(\frac{m\omega_m^q}{n\omega_n}\right)^{
		\frac{n}{n-m}}
	\left(\frac{V(K)}{\omega_n}\right)^{
		\frac{m(q-1)}{n-m}}.
	\]
	On the other hand, 
	\[
	V(\mathbb I_{m,q}B^n)
	=
	\omega_n
	\left(\frac{m\omega_m^q}{n\omega_n}\right)^{
		\frac{n}{n-m}}.
	\]
	Combining these identities yields the desired inequality.
	
	Since $1<q<n$ and $L_0$ is a ball centered at the origin,
	the equality conditions in Theorem~\ref{thm-6-7} and Theorem~\ref{thm-6-8} imply that equality
	holds if and only if $K$ is a dilate of $L_0$, that is,
	a Euclidean ball centered at the origin.
\end{proof}

\noindent\textbf{Remark 6.10.}
The range $1\le q\le n$ in Theorem~\ref{thm-6-7} is sharp,
even for origin-centered ellipsoids.

First suppose that $q<1$. Let $L=B^n$, and choose an
origin-centered ellipsoid $K$ that is not a ball and satisfies
$V(K)=\omega_n$. Since $K$ is origin-symmetric,
$s_K(u)=2\rho_K(u)$, and hence
\[
\frac{1}{n\omega_n}\int_{S^{n-1}}s_K(u)^n\,du
=
\frac{2^n}{n\omega_n}\int_{S^{n-1}}\rho_K(u)^n\,du
=
2^n.
\]
The function $t\mapsto t^{(q-1)/(n-1)}$ is strictly convex on
$(0,\infty)$. Applying Jensen's inequality to $s_K^n$ with
respect to the probability measure $du/(n\omega_n)$ gives
\[
\begin{aligned}
	\frac{1}{n\omega_n}
	\int_{S^{n-1}}s_K(u)^{\frac{n(q-1)}{n-1}}\,du
	>
	\left(
	\frac{1}{n\omega_n}\int_{S^{n-1}}s_K(u)^n\,du
	\right)^{\frac{q-1}{n-1}}=2^{\frac{n(q-1)}{n-1}},
\end{aligned}
\]
where the inequality is strict because $s_K$ is continuous
and nonconstant. Since $s_{B^n}(u)=2$, it follows that
\[
\int_{S^{n-1}}s_K(u)^{\frac{n(q-1)}{n-1}}\,du
>
\int_{S^{n-1}}s_{B^n}(u)^{\frac{n(q-1)}{n-1}}\,du.
\]
Using the volume formula for $\mathbb I_{1,q}$, we obtain
\[ 
	V(\mathbb I_{1,q}(K,B^n))>
	V(\mathbb I_{1,q}(B^n,B^n))=
	\left(\frac{2^n}{n\omega_n}\right)^{\frac n{n-1}}\omega_n.
\]
Because $V(K)=V(B^n)=\omega_n$, the last expression is exactly
the right-hand side of the inequality in
Theorem~\ref{thm-6-7} for the pair $(K,B^n)$.
Thus, that inequality does not extend to $q<1$.

For $q>n$, take $K=B^n$ and choose an origin-centered
ellipsoid $L$ that is not a ball and satisfies
$V(L)=\omega_n$. Applying the same argument to $s_L^n$
with the strictly convex function
$t\mapsto t^{(n-q)/(n-1)}$ yields
\[
\frac{1}{n\omega_n}
\int_{S^{n-1}}s_L(u)^{\frac{n(n-q)}{n-1}}\,du
>
2^{\frac{n(n-q)}{n-1}}.
\]
Consequently,
\[
V(\mathbb I_{1,q}(B^n,L))
>
V(\mathbb I_{1,q}(B^n,B^n)).
\]
Since $V(L)=V(B^n)$, this shows that the inequality in
Theorem~\ref{thm-6-7} does not extend to $q>n$.\\

	\textbf{Proof of Theorem~\ref{thm-1-3}.} It follows immediately from Theorem~\ref{thm-6-7} and Theorem~\ref{thm-6-8}.\quad $\square$
	
	\medskip
	\noindent\textbf{Acknowledgements.}
	This work was supported by the National Natural Science Foundation of China (Grant No.~12322103).
    
    \medskip
	\noindent\textbf{Conflict of Interest Statement.}
	The author declares that there are no conflicts of interest.
	
	\medskip
	\noindent\textbf{Data Availability Statement.}
	No data was used for the research described in the article.

\bigskip
\bigskip

\noindent
{\scshape Department of Mathematics, Shanghai University, Shanghai 200444, China}

\medskip

\noindent
\textit{Email address:}
\href{mailto:yrenqing2020@163.com}{yrenqing2020@163.com}

\end{document}